\AddToHook{begindocument/before}{} 
\DeclareHookRule{begindocument}{hyperref}{before}{my-mdpi} 

\documentclass[preprints,article,accept,moreauthors,pdftex]{my-mdpi}

\usepackage[normalem]{ulem} 
\usepackage{amsfonts} 
\usepackage{wasysym} 
\usepackage{mathcomp} 
\usepackage{pifont} 
\usepackage{bm} 

\firstpage{1}
\articlenumber{2503}
\doinum{10.3390/math11112503} 
\pubvolume{11}
\issuenum{11}
\pubyear{2023}
\copyrightyear{2023}
\externaleditor{Academic Editors: Juan Eduardo N{\fontencoding{T5}\selectfont{\'a}}poles Valdes and Miguel Vivas-Cortez}
\datereceived{20 April 2023}
\daterevised{25 May 2023}
\dateaccepted{25 May 2023}
\datepublished{29 May 2023}
\hreflink{https://\linebreak doi.org/10.3390/math11112503}
\usepackage[T5,T1]{fontenc}
\usepackage[english]{babel}

\pdfoutput=1

\Title{Numerical Investigation of the Fractional Oscillation Equations under the Context of Variable Order Caputo Fractional Derivative via Fractional Order Bernstein Wavelets}

\TitleCitation{Numerical Investigation of the Fractional Oscillation Equations under the Context of Variable Order Caputo Fractional Derivative via Fractional Order Bernstein Wavelets}

\Author{Ashish~Rayal~\textsuperscript{1}\href{https://orcid.org/0000-0001-8954-6605}{\orcidicon}, 
Bhagawati Prasad~Joshi~\textsuperscript{2}\href{https://orcid.org/0000-0001-6935-5975}{\orcidicon}, 
Mukesh~Pandey~\textsuperscript{3} 
and Delfim F. M.~Torres~\textsuperscript{4}\textsuperscript{,}*\href{https://orcid.org/0000-0001-8641-2505}{\orcidicon}}

\AuthorNames{Ashish Rayal, Bhagawati Prasad Joshi, Mukesh Pandey and Delfim F. M. Torres}

\AuthorCitation{Rayal, A.; Joshi, B.P.; Pandey, M.; Torres, D.F.M.}

\address{\textsuperscript{1} \quad Department of Mathematics, 
School of Applied and Life Sciences, Uttaranchal University, 
Dehradun~248007,~India; 
ashishrayal@uttaranchaluniversity.ac.in

\textsuperscript{2} \quad Department of Mathematics, 
Graphic Era Hill University, Bhimtal 263136, India; 
bpjoshi@gehu.ac.in

\textsuperscript{3} \quad School of Computer Science, 
University of Petroleum \& Energy Studies, Dehradun~248007,~India; 
mukesh.110759@stu.upes.ac.in

\textsuperscript{4} \quad Center for Research and Development in Mathematics and Applications~(CIDMA), 
Department~of~Mathematics, University of Aveiro, 3810-193 Aveiro, Portugal}

\corres{Correspondence: delfim@ua.pt}

\abstract{This article describes an approximation technique based on fractional order Bernstein 
wavelets for the numerical simulations of fractional oscillation equations under variable order, 
and the fractional order Bernstein wavelets are derived by means of fractional Bernstein polynomials. 
The oscillation equation describes electrical circuits and exhibits a wide range of nonlinear dynamical behaviors. 
The proposed variable order model is of current interest in a lot of application areas in engineering and applied sciences. 
The purpose of this study is to analyze the behavior of the fractional force-free and forced oscillation equations 
under the variable-order fractional operator. The basic idea behind using the approximation technique 
is that it converts the proposed model into non-linear algebraic equations with the help of collocation 
nodes for easy computation. Different cases of the proposed model are examined under the selected variable 
order parameters for the first time in order to show the precision and performance of the mentioned scheme. 
The dynamic behavior and results are presented via tables and graphs to ensure the validity of the mentioned scheme. 
Further, the behavior of the obtained solutions for the variable order is also depicted. From the calculated results, 
it is observed that the mentioned scheme is extremely simple and efficient for examining the behavior 
of nonlinear random (constant or variable) order fractional models occurring in engineering and science.}



\begin{document}

\makeatletter

	
\section{Introduction}
\label{sect:sec1-mathematics-2384216}

In previous decades, concepts of fractional order calculus (FOC) have been extensively 
employed in all areas of science, economics, and engineering fields, and they are growing 
very fast in developing and describing the behavior of models due to their relation to hereditary, fractals, and memory~\cite{B1-mathematics-2384216,B2-mathematics-2384216,B3-mathematics-2384216,B4-mathematics-2384216}. 
FOC also gives several fractional-order integral and derivative operators and numerical solutions with high accuracy. 
The classifications of fractional operators are based on the concepts of the singular kernel, non-singular kernel, 
nonlocal kernel, and non-singular kernel. Some of them are Caputo, Atangana-Baleanu, Caputo-Fabrizio, Riesz, 
Riemann-Liouville, and Hadamard. For example, the authors in~\cite{B5-mathematics-2384216} introduced 
the operational matrices of fractional Bernstein functions to solve fractional differential equations (FDEs), 
and Alshbool et~al.~\cite{B6-mathematics-2384216} proposed the concept of operational matrices based 
on fractional Bernstein functions for solving integro-differential equations under the Caputo operator. 
The use of new fractional operators in the geometry of real-world models has made significant advancements 
in this domain~\cite{B7-mathematics-2384216,B8-mathematics-2384216}. In most cases, the researchers have not 
achieved desirable solutions using integer-order operators. This fact emphasizes the significance 
of new differential operators in modeling real-world problems.

The most extended area for FOC involves variable-ordered operators because the order of fractional 
operators could be any arbitrary value. The fractional operators under variable order override 
the phenomenon of constant-order fractional operators. This encourages us to investigate 
some new concepts in the proposed manner due to their numerous application areas in engineering and science. 
The nonlocal characteristics of systems are more apparent with non-constant-order fractional calculus. 
The FOC with variable order is used to model many phenomena such as anomalous diffusions with constant 
and variable orders, viscoelastic spherical indentation, transient dispersion in heterogeneous media, alcoholism, and so on~\cite{B9-mathematics-2384216,B10-mathematics-2384216,B11-mathematics-2384216,B12-mathematics-2384216}. 
It is usually more complicated to estimate the explicit solution of fractional differential equations (FDEs) 
under variable order. Hence, it is necessary to describe numerical approaches for the solution of such problems. 
There are several schemes for solving FDEs in variable order. Among these schemes, wavelet-related schemes 
are more attractive and efficient for solving this type of problem due to wavelets’ important features 
like compact support, spectral accuracy, orthogonality, and localization.

Wavelets~\cite{B13-mathematics-2384216,B14-mathematics-2384216} are the good localized and oscillatory 
functions that give the basis for several spaces. In approximation theory, there is lots of literature 
available concerning the power series and Fourier series. The approximation of an arbitrary function 
through wavelet polynomials is a recent development in approximation theory. The wavelet expansion 
is more generalized than any other expansion, such as the power and Fourier series. The main reason 
for the discovery of wavelets is that the Fourier series cannot analyze the signal in both the 
frequency and time domains. The important benefit of the wavelet transform is its ability to analyze 
the signal simultaneously in the frequency and time domains. Orthogonal wavelets play an important 
role in solving differential and integral equations. In the past two decades, wavelet approaches 
have been extensively employed to solve differential equations of arbitrary order arising in numerous 
engineering and scientific problems. Several researchers have used wavelet-based approximation 
approaches to solve different classes of differential equations. 
See these references~\cite{B15-mathematics-2384216,B16-mathematics-2384216,B17-mathematics-2384216,%
B18-mathematics-2384216,B19-mathematics-2384216,B20-mathematics-2384216} for more applications of wavelets.

Here, we introduce the application of FOC under variable-order 
for the modeling of nonlinear oscillation equations as
\begin{equation}
\label{eq:FD1-mathematics-2384216}
\text{D}_{0,\text{t}}^{\text{$\upalpha$}(\text{t})}\Im(\text{t}) - \text{$\upmu$}\Im^{\prime}(\text{t}) 
+ \text{$\upmu$}\Im^{\prime}(\text{t})\Im^{2}(\text{t}) + \text{a}\Im(\text{t}) 
+ \text{b}\Im^{3}(\text{t}) = \Phi\left( {\text{$\upomega$},\text{f},\text{t}} \right);
\ \text{$\upalpha$}(\text{t}) \in \left( {1,2} \right\rbrack,
\tag{1}
\end{equation}
with the initial value conditions
\begin{equation}
\nonumber\label{eq:FD2-mathematics-2384216}
\Im(0) = 1,\ \Im^{\prime}(0) = 0,
\end{equation}
where $\Phi\left( {\text{f},\text{$\upomega$},\text{t}} \right) $ is the forcing term 
or prescribed excitation, $\text{$\upomega$} $ is the driving force’s angular frequency, 
$\Im(\text{t}) $ is the system response, $f$ is the amplitude of the excitation, $a$ \& $b$ 
are constant parameters, $\text{$\upmu$} $ is the damping parameter of the considered system, 
and $\text{D}_{0,\text{t}}^{\text{$\upalpha$}(\text{t})} $ is a fractional 
Caputo derivative with order $\text{$\upalpha$}(\text{t})$.

The primary aim of the present work is to estimate a more convenient wavelet solution of the 
fractional oscillation equation under variable order via the fractional order Bernstein wavelets (FOBWs) basis. 
The proposed method involves approximating the unknown function using a truncated FOBWs basis. 
After approximating this function, a series of nonlinear algebraic equations 
is formed for estimating the wavelet coefficient vector.

This work is significantly helpful for the study of any type of variable-order nonlinear fractional model. 
Some of the advantages of this work are listed as follows:
\begin{enumerate}[label=$\bullet$]
\item The present scheme works for the first time with the Caputo fractional 
derivative under variable order in the introduced model. This work deals with 
the replacement of constant order by variable order in the considered 
nonlinear model under the fractional operator.

\item From a computational point of view, only fewer terms of FOBWs bases are applied 
to achieve very satisfactory and effective results in comparison to existing methods, 
which is a key feature of the mentioned scheme.

\item The introduced FOBWs are simple bases from a computational point of view; therefore, 
these bases could be seen as a convenient and appropriate tool in this work 
for solving the fractional oscillation equation under variable order.

\item The mentioned scheme is very easy to implement and provides 
better accuracy in comparison to other existing schemes.

\item The present study is very useful to investigate the behavior 
of several nonlinear variable-order fractional models with fewer errors.
\end{enumerate}

The remaining portion of the manuscript is designed as follows: Section~\ref{sect:sec2-mathematics-2384216} 
provides the basic preliminaries about fractional operators and special functions. Section~\ref{sect:sec3-mathematics-2384216} 
recalls the related work. The definition of FOBW’s basis is given in Section~\ref{sect:sec4-mathematics-2384216}. 
In Section~\ref{sect:sec5-mathematics-2384216}, the approximation of function through FOBWs has been explained. 
Section~\ref{sect:sec6-mathematics-2384216} presents the FOBWs scheme for the evaluation of the fractional oscillation 
equation under variable order. Section~\ref{sect:sec7-mathematics-2384216} shows the result of the convergence analysis. 
In Section~\ref{sect:sec8-mathematics-2384216}, some applications on different parameters are evaluated, 
which illustrates the efficiency of the mentioned approach. 
The conclusion is drawn in Section~\ref{sect:sec9-mathematics-2384216}.


\section{Preliminaries}
\label{sect:sec2-mathematics-2384216}

In this study, the following concepts of variable order fractional operators and special functions are used.

\vspace{12pt}
\noindent \textbf{\boldmath{Definition}} \textbf{\boldmath{1.}} 
\emph{The fractional Caputo differentiation of} $\Im(\text{t}) \in \text{L}^{2}\left\lbrack {0,1} 
\right\rbrack $ \emph{with order} $\text{$\upalpha$}(\text{t}) $ 
\emph{is given~by~\cite{B21-mathematics-2384216}.}
\begin{equation}
\label{eq:FD3-mathematics-2384216}
\text{D}_{0,\text{t}}^{\text{$\upalpha$}(\text{t})}\Im(\text{t}) 
= 
\begin{cases}
{\frac{1}{(1 + \text{n} - \text{$\upalpha$}(\text{t})!)}{
\int\limits_{0}^{\text{t}}{\Im^{(\text{n})}(\text{$\uptau$}){(\text{t} 
- \text{$\uptau$})}^{\text{n} - \text{$\upalpha$}(\text{t}) - 1}}}\text{d}\text{$\uptau$},} 
& {\text{n} - 1 < \ \text{$\upalpha$}(\text{t}) < \text{n} \in \mathbb{N}} \\
{\Im^{(\text{n})}(\text{t}),} & ~~{\text{$\upalpha$}(\text{t}) = \text{n} \in \mathbb{N}} \\
\end{cases}.
\tag{2}
\end{equation}

\vspace{12pt}
\noindent \textbf{\boldmath{Definition}} \textbf{\boldmath{2.}} 
\emph{The fractional Riemann-Liouville integral of} $\Im(\text{t}) $\emph{with order} 
$\text{$\upalpha$}(\text{t}) $ \emph{is given as~\cite{B21-mathematics-2384216}.}
\begin{equation}
\label{eq:FD4-mathematics-2384216}
\text{I}_{0,\text{t}}^{\text{$\upalpha$}(\text{t})}\Im(\text{t}) 
= \frac{1}{(1 + \text{$\upalpha$}(\text{t})!)}{\int\limits_{0}^{\text{t}}{(\text{t} 
- \text{$\uptau$})}^{\text{$\upalpha$}(\text{t}) - 1}}
\Im(\text{$\uptau$})\text{d}\text{$\uptau$},\quad 0 < \text{t}.
\tag{3}
\end{equation}

\vspace{12pt}
In addition, the connections between fractional Caputo derivatives and fractional 
Riemann-Liouville integral for $\text{n} - 1 < \ \text{$\upalpha$}(\text{t}) \leq \text{n} $ 
and $\text{$\uplambda$} > 0 $ are:
\begin{equation}
\label{eq:FD5-mathematics-2384216}
\text{I}_{0,\text{t}}^{\text{$\uplambda$}}(\text{D}_{0,\text{t}}^{\text{$\uplambda$}}\Im(\text{t})) 
= \Im(\text{t}) - {\sum\limits_{\text{j} = 0}^{{\lceil\text{$\uplambda$}\rceil} - 1}
\frac{\text{t}^{\text{j}}}{\text{j}!}}\Im^{(\text{j})}(0),\quad 0 < \text{t}.
\tag{4}
\end{equation}
\begin{equation}
\label{eq:FD6-mathematics-2384216}
\text{I}_{0,\text{t}}^{\text{n} - \text{$\upalpha$}(\text{t})}(\Im^{(\text{n})}(\text{t})) 
= \text{D}_{0,\text{t}}^{\text{$\upalpha$}(\text{t})}\Im(\text{t}) - {\sum\limits_{\text{j} 
= {\lceil{\text{$\upalpha$}(\text{t})}\rceil}}^{\text{n} - 1}\frac{\text{t}^{\text{j} 
- \text{$\upalpha$}(\text{t})}}{(\text{j} 
- \text{$\upalpha$}(\text{t}))!}}\Im^{(\text{j})}(0),\quad 0 < \text{t}.
\tag{5}
\end{equation}

Further, the relationship between factorial, gamma function, 
and binomial coefficients is given as~\cite{B22-mathematics-2384216}

\vspace{-3pt}
\begin{equation}
\nonumber\label{eq:FD7-mathematics-2384216}
(\text{j} - \text{$\upalpha$}(\text{t}))! 
= \text{$\Gamma$}\left( {1 + \text{j} - \text{$\upalpha$}(\text{t})} \right).
\end{equation}
\begin{equation}
\nonumber\label{eq:FD8-mathematics-2384216}
\begin{pmatrix}
\text{a}_{1} \\
\text{a}_{2} \\
\end{pmatrix} = \frac{(\text{a}_{1})!}{(\text{a}_{2})!(\text{a}_{1} 
- \text{a}_{2})!},\quad\text{a}_{1} \geq \text{a}_{2}.
\end{equation}

Here, $\text{$\Gamma$}\left( \cdot \right) $ denotes the gamma function, and $(\cdot)!$ 
denotes the factorial function, and $\left\lceil \text{$\uplambda$} \right\rceil $ 
is the well-known ceiling function or least integer function. For the proposed model, 
we take $n = 2$ in the above definitions, so that $\text{$\upalpha$}(\text{t}) 
\in \left( {1,2} \right\rbrack$.


\section{Related Work}
\label{sect:sec3-mathematics-2384216}

The oscillation equation is the most classical differential equation in nonlinear dynamics that models systems under self-sustained oscillation and is used as a model in image processing, neurology, electronics, and so on~\cite{B23-mathematics-2384216,B24-mathematics-2384216,B25-mathematics-2384216}. Various numerical and analytical approaches have been introduced for solving oscillation equations. Cordshooli and Vahidi~\cite{B26-mathematics-2384216} proposed the series solution of the oscillation equation by using the adomian decomposition scheme (ADS). In~\cite{B27-mathematics-2384216}, Vahidi et~al. employed restarted ADS to solve the oscillation equation. In~\cite{B28-mathematics-2384216}, Doha et~al. presented a collocation scheme combined with an ultraspherical wavelet for approximating the oscillation equation. In~\cite{B29-mathematics-2384216}, the authors presented an efficient solution of the fractional oscillation equation through a modified Legendre wavelet. Khan~\cite{B30-mathematics-2384216} presented the approximate solution of the oscillation equation through the homotopy perturbation method. In~\cite{B31-mathematics-2384216}, Kumar and Varshney proposed the numerical simulation of the Vander Pol equation through the Lindstedt-Poincare scheme. Recently, Hamed et~al.~\cite{B32-mathematics-2384216} provided a numerical treatment of the stochastic oscillation equation using the Wiener--Hermite expansion approach.

Many physical and biological problems are governed through FDEs under variable order, such as the cable equation~\cite{B33-mathematics-2384216}, the Rayleigh--Stokes equation~\cite{B34-mathematics-2384216}, the Schrödinger equation~\cite{B35-mathematics-2384216}, and so on. The explicit solutions to most of the FDEs in variable order are difficult to find. Therefore, obtaining solutions to such problems has taken the attention of several researchers. A detailed summary of the solutions of FDEs under variable order arising in the fields of biology, engineering, and physics is given in Table~\ref{tabref:mathematics-2384216-t001}. It has been revealed from the literature review that analysis of the mathematical, engineering, and physical models associated with variable-order fractional derivatives rather than derivatives of integer order provides highly significant results.    

The oscillation equation has only been solved for fractional constant order, but in this paper, we introduce the oscillation equation under the concept of variable order fractional derivative due to the advantages of employing variable fractional order. In order to more efficiently solve the fractional oscillation equation under variable order, the FOBWs are introduced in this study. The present study aims to extend the applications of FOBWs with collocation techniques to the approximate solutions of fractional oscillation equations under variable order and analyze their behavior with different parameters. The computing complexity of the algebraic set can be decreased due to the structural redundancy of the FOBWs. The errors under several fractional variable orders are computed, which proves the effectiveness of the scheme mentioned. So, keeping all the facts in mind and influenced by the good performance of the above-mentioned approaches, we will employ an effective wavelet approach for the numerical analysis of a variable-order nonlinear fractional model such as the fractional oscillation equation.

\begin{table}[H]
\tablesize{\small}
\caption{Detail of the numerical schemes for the solutions of FDEs under variable order.}
\label{tabref:mathematics-2384216-t001}
\begin{adjustwidth}{-\extralength}{0cm}
\setlength{\cellWidtha}{\fulllength/3-2\tabcolsep-1.8in}
\setlength{\cellWidthb}{\fulllength/3-2\tabcolsep-0in}
\setlength{\cellWidthc}{\fulllength/3-2\tabcolsep+1.8in}
\scalebox{1}[1]{\begin{tabularx}{\fulllength}{>{\centering\arraybackslash}m{\cellWidtha}>{\centering\arraybackslash}m{\cellWidthb}>{\centering\arraybackslash}m{\cellWidthc}}
\toprule
\textbf{S. No.} & \textbf{Authors/References} & \textbf{Description of Schemes}\\
\cmidrule{1-3}
1 & Xu and Erturk {\cite{B36-mathematics-2384216}} & Proposed a finite difference approach for the solution of an integro-differential equation under variable fractional order.\\
\cmidrule{1-3}
2 & Wang and Vong {\cite{B37-mathematics-2384216}} & Presented the difference approach for the reformed version of the anomalous fractional wave equation.\\
\cmidrule{1-3}
3 & Fu, Chen and Ling {\cite{B38-mathematics-2384216}} & Introduced a novel numerical scheme for approximate solutions to random order fractional diffusion models.\\
\cmidrule{1-3}
4 & Zayernouri and Karniadakis {\cite{B39-mathematics-2384216}} & Described the fractional spectral collocation technique for nonlinear partial differential equations under variable fractional order.\\
\cmidrule{1-3}
5 & Chen, Wei, Liu and Yu {\cite{B40-mathematics-2384216}} & Developed the Legendre wavelet technique for the solution of nonlinear random-order FDEs.\\
\cmidrule{1-3}
6 & Yaghoobi, Moghaddam, and Ivaz {\cite{B41-mathematics-2384216}} & Provide an efficient solution for variable-order time delay FDEs through cubic spline approximation.\\
\cmidrule{1-3}
7 & Aguilar, Coronel-Escamilla, Gomez-Aguilar, Alvarado-Martinez, and Romero-Ugalde {\cite{B42-mathematics-2384216}} & Presented the simulation of FDEs under variable order with the Mittag-Leffler kernel by the artificial neural network technique.\\
\cmidrule{1-3}
8 & Heydari {\cite{B43-mathematics-2384216}} & Proposed an approach for a fractional variable-order optimal control model in the Atangana-Baleanu sense with the help of the Chebyshev cardinal functions.\\
\cmidrule{1-3}
9 & Nemati, Lima, and Torres {\cite{B44-mathematics-2384216}} & Introduce an approach for FDEs under variable order through Bernoulli~polynomials.\\
\cmidrule{1-3}
10 & Kaabar, Refice, Souid, Mart{\fontencoding{T5}\selectfont{\'i}}nez, Etemad, Siri, and Rezapour {\cite{B45-mathematics-2384216}} & Established the stability criteria for the solution of the fractional boundary problem under variable order.\\
\bottomrule
\end{tabularx}}
\end{adjustwidth}
\end{table}


\section{Development of Fractional Order Bernstein Wavelets}
\label{sect:sec4-mathematics-2384216}

In the current section, first the definition of fractional order Bernstein polynomials is recalled, and then the Bernstein wavelets are constructed in fractional form.


\subsection{Fractional Order Bernstein Polynomials}
\label{sect:sec4dot1-mathematics-2384216}

\textls[-15]{The fractional order Bernstein polynomials of order $\text{$\upupsilon$}\text{$\upgamma$} $ 
are defined in explicit form~as~\cite{B46-mathematics-2384216}:}

\vspace{-6pt}
\begin{equation}
\label{eq:FD9-mathematics-2384216}
\text{B}_{\text{$\upupsilon$},\text{M}}^{\text{$\upgamma$}}(\text{t}) = \sqrt{1 + 2\text{M} - 2\text{$\upupsilon$}}\quad\left( {1 - \text{t}^{\text{$\upgamma$}}} \right)^{\text{M} - \text{$\upupsilon$}}{\sum\limits_{\text{i} = 0}^{\text{$\upupsilon$}}\left( {- 1} \right)^{\text{i}}}\begin{pmatrix}
{1 + 2\text{M} - \text{i}} \\
{\text{$\upupsilon$} - \text{i}} \\
\end{pmatrix}\begin{pmatrix}
\text{$\upupsilon$} \\
\text{i} \\
\end{pmatrix}\text{t}^{\text{$\upgamma$}(\text{$\upupsilon$} - \text{i})}.
\tag{6}
\end{equation}

The above polynomials in Equation (6) are orthogonal under the weighted function $\text{$\Omega$}(\text{t}) = \text{t}^{\text{$\upgamma$} - 1} $ on [0, 1] as
\begin{equation}
\label{eq:FD10-mathematics-2384216}
{\int\limits_{0}^{1}{\text{B}_{\text{$\upupsilon$},\text{M}}^{\text{$\upgamma$}}(\text{t})\text{B}_{\text{$\upvartheta$},\text{M}}^{\text{$\upgamma$}}(\text{t})}}\text{$\Omega$}(\text{t})\text{d}\text{t} = \begin{cases}
{0,} & {\text{$\upupsilon$} \neq \text{$\upvartheta$}} \\
{1/\text{$\upgamma$},} & {\text{$\upupsilon$} = \text{$\upvartheta$}} \\
\end{cases}.
\tag{7}
\end{equation}

In addition, the other form of the above polynomials is given as
\begin{equation}
\label{eq:FD11-mathematics-2384216}
\text{B}_{\text{$\upupsilon$},\text{M}}^{\text{$\upgamma$}}(\text{t}) = \sqrt{1 + 2\text{M} - 2\text{$\upupsilon$}}\ {\sum\limits_{\text{i} = 0}^{\text{$\upupsilon$}}\left( {- 1} \right)^{\text{i}}}\frac{\begin{pmatrix}
{1 + 2\text{M} - \text{i}} \\
{\text{$\upupsilon$} - \text{i}} \\
\end{pmatrix}\begin{pmatrix}
\text{$\upupsilon$} \\
\text{i} \\
\end{pmatrix}}{\begin{pmatrix}
{\text{M} - \text{i}} \\
{\text{$\upupsilon$} - \text{i}} \\
\end{pmatrix}}{\widetilde{\text{B}}}_{\text{$\upupsilon$} - \text{i},\text{M} - \text{i}}^{\text{$\upgamma$}}(\text{t}),
\tag{8}
\end{equation}
where
\begin{equation}
\nonumber\label{eq:FD12-mathematics-2384216}
{\widetilde{\text{B}}}_{\text{$\upupsilon$},\text{M}}^{\text{$\upgamma$}}(\text{t}) = {\sum\limits_{\text{i} = 0}^{\text{M} - \text{$\upupsilon$}}\left( {- 1} \right)^{\text{i}}}\begin{pmatrix}
\text{M} \\
\text{$\upupsilon$} \\
\end{pmatrix}\begin{pmatrix}
{\text{M} - \text{$\upupsilon$}} \\
\text{i} \\
\end{pmatrix}\text{t}^{\text{$\upgamma$}(\text{$\upupsilon$} + \text{i})}.
\end{equation}

In Equation (8), '$i$' is a whole number that represents the index value of the given summation, and $M$ is a natural number.


\subsection{Fractional Order Bernstein Wavelets}
\label{sect:sec4dot2-mathematics-2384216}

The FOBWs $\Xi_{\text{$\upeta$},\text{$\upupsilon$}}^{\text{$\upgamma$}}(\text{t}) = \Xi(\text{k},\text{t},\text{$\upupsilon$},\text{$\upeta$},\text{$\upgamma$}) $ have the arguments: k is a natural number, t represents time, $\text{$\upupsilon$} $ is the order of Bernstein polynomial such that $\text{$\upupsilon$} = 0,1,2,3,\ldots,\text{M} \in \mathbb{N} $, $\text{$\upeta$} = 1,2,3,\ldots,2^{\text{k} - 1}, $ and $\text{$\upgamma$} > 0. $

The FOBWs is defined on $\left\lbrack {0,1} \right\rbrack $ as
\begin{equation}
\label{eq:FD13-mathematics-2384216}
\Xi_{\text{$\upeta$},\text{$\upupsilon$}}^{\text{$\upgamma$}}(\text{t}) = \begin{cases}
{\sqrt{\text{$\upgamma$}}2^{\frac{\text{k} - 1}{2}}\text{B}_{\text{$\upupsilon$},\text{M}}^{\text{$\upgamma$}}(1 + 2^{\text{k} - 1}\text{t} - \text{$\upeta$}),} & {\frac{\text{$\upeta$} - 1}{2^{\text{k} - 1}} \leq \text{t} \leq \frac{\text{$\upeta$}}{2^{\text{k} - 1}}} \\
{0,} & {\text{o}\text{t}\text{h}\text{e}\text{r}\text{w}\text{i}\text{s}\text{e}} \\
\end{cases},
\tag{9}
\end{equation}
where $\text{B}_{\text{$\upupsilon$},\text{M}}^{\text{$\upgamma$}}(\text{t}) $ is the fractional order Bernstein polynomials of order $\text{$\upupsilon$}\text{$\upgamma$} $ defined in Section~\ref{sect:sec4dot1-mathematics-2384216}.
\begin{enumerate}[label=$\bullet$]
\item The set of FOBWs forms the orthonormal set on $\left\lbrack {0,1} \right\rbrack $ under the weighted function $\text{$\Omega$}_{\text{k},\text{$\upeta$}}(\text{t}), $ where

\vspace{-3pt}
\begin{equation}
\label{eq:FD14-mathematics-2384216}
\text{$\Omega$}_{\text{k},\text{$\upeta$}}(\text{t}) = \text{$\Omega$}(1 + 2^{\text{k} - 1}\text{t} - \text{$\upeta$}).
\tag{10}
\end{equation}
\end{enumerate}
\noindent i.e.,
\begin{equation}
\nonumber\label{eq:FD15-mathematics-2384216}
{\int\limits_{0}^{1}{\Xi_{\text{$\upeta$},\text{$\upupsilon$}}^{\text{$\upgamma$}}(\text{t})\Xi_{\text{$\upeta$},\text{$\upvartheta$}}^{\text{$\upgamma$}}(\text{t})}}\text{$\Omega$}_{\text{k},\text{$\upeta$}}(\text{t})\text{d}\text{t} = \begin{cases}
{0,} & {\text{$\upupsilon$} \neq \text{$\upvartheta$}} \\
{1,} & {\text{$\upupsilon$} = \text{$\upvartheta$}} \\
\end{cases}.
\end{equation}

\begin{enumerate}[label=$\bullet$]
\item The FOBWs have compact support, i.e.,\begin{equation}
\nonumber\label{eq:FD16-mathematics-2384216}
\begin{array}{r}
{\text{s}\text{u}\text{p}\text{p}\left( {\Xi_{\text{$\upeta$},\text{$\upupsilon$}}^{\text{$\upgamma$}}(\text{t})} \right) = \left\{ \overline{\text{t}:\Xi_{\text{$\upeta$},\text{$\upupsilon$}}^{\text{$\upgamma$}}(\text{t})\neq 0} \right\}}\vspace{4pt}\\
{= \left\lbrack {\frac{\text{$\upeta$} - 1}{2^{\text{k} - 1}},\frac{\text{$\upeta$}}{2^{\text{k} - 1}}} \right\rbrack.} \\
\end{array}
\end{equation}

\item The FOBWs basis is exactly the classical Bernstein wavelets for unit $\text{$\upgamma$}. $
\end{enumerate}

\begin{figure}[H]
\includegraphics[scale=1.4]{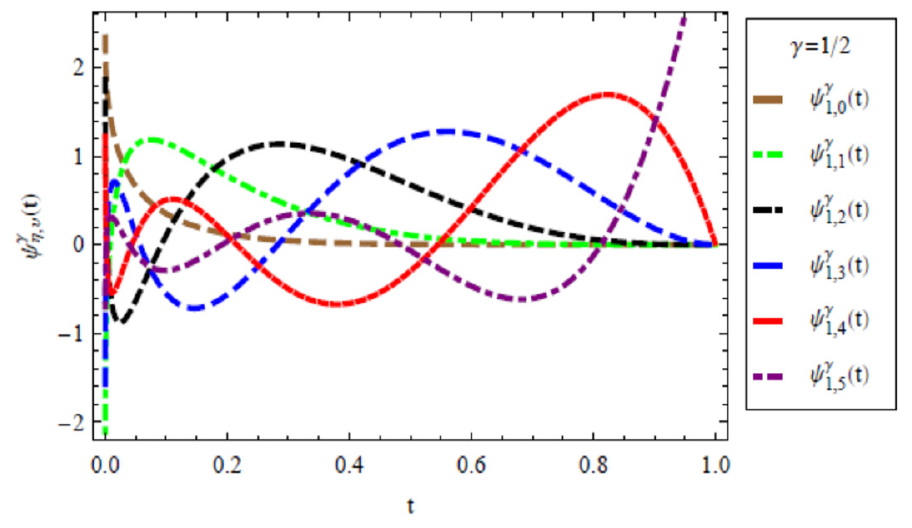}
\caption{Plot of the FOBWs for M = 5, k = 1, and $\text{$\upgamma$} = 1/2. $}
\label{fig:mathematics-2384216-f001}
\end{figure}

The FOBWs are displayed in Figure~\ref{fig:mathematics-2384216-f001} for k = 1, M = 5, and $\text{$\upgamma$} = 1/2. $    

Now, using the above wavelet basis, the approximation of any function in the Hilbert space is stated in next section.


\section{Function Approximations via Introduced FOBWs}
\label{sect:sec5-mathematics-2384216}

In order to employ FOBWs for solving the proposed model, we need to map the unknown function to FOBWs. For this, the concept of function approximation is used. The unknown function can be approximated by this concept in terms of a known wavelet function with wavelet coefficients.

Any arbitrary function $\Im(\text{t}) \in \text{L}^{2}\left\lbrack {0,1} \right\rbrack $ can be formulated in a combination of \linebreak FOBWs~as~\cite{B13-mathematics-2384216,B40-mathematics-2384216,B47-mathematics-2384216}

\vspace{-4pt}
\begin{equation}
\label{eq:FD17-mathematics-2384216}
\Im(\text{t}) \approx {\sum\limits_{\text{$\upeta$} = 1}^{\infty}{\sum\limits_{\text{$\upupsilon$} = 0}^{\infty}{\text{u}_{\text{$\upeta$},\text{$\upupsilon$}}\Xi_{\text{$\upeta$},\text{$\upupsilon$}}^{\text{$\upgamma$}}(\text{t})}}}.
\tag{11}
\end{equation}

For approximation purposes, the truncated form of Equation (11) is written as
\begin{equation}
\label{eq:FD18-mathematics-2384216}
\Im(\text{t}) \approx {\sum\limits_{\text{$\upeta$} = 1}^{2^{\text{k} - 1}}{\sum\limits_{\text{$\upupsilon$} = 0}^{\text{M}}{\text{u}_{\text{$\upeta$},\text{$\upupsilon$}}\Xi_{\text{$\upeta$},\text{$\upupsilon$}}^{\text{$\upgamma$}}(\text{t}) = \text{U}^{\text{T}}}}}\text{$\Psi$}_{\widetilde{\text{$\upsigma$}}}(\text{t}),
\tag{12}
\end{equation}
where $\text{U} $ be the unknown wavelet coefficients associated with FOBWs $\text{$\Psi$}_{\widetilde{\text{$\upsigma$}}}^{\text{$\upgamma$}}(\text{t}) $ given by
\begin{equation}
\label{eq:FD19-mathematics-2384216}
\text{U}^{\text{T}} = \ \left\lbrack 
\begin{array}{l}
{\text{u}_{1,0},\text{u}_{1,1},\text{u}_{1,2},\ldots,\text{u}_{1,\text{M} - 1},\text{u}_{2,0},\text{u}_{2,1},\text{u}_{2,2},\ldots,\text{u}_{2,\text{M} - 1},\ldots,} \vspace{3pt}\\
{\text{u}_{2^{\text{k} - 1},0},\text{u}_{2^{\text{k} - 1},1},\ldots,\text{u}_{2^{\text{k} - 1},\text{M}}} \\
\end{array} \right\rbrack,
\tag{13}
\end{equation}
\begin{equation}
\label{eq:FD20-mathematics-2384216}
\text{$\Psi$}_{\widetilde{\text{$\upsigma$}}}^{\text{$\upgamma$}}(\text{t}) = \left\lbrack \begin{array}{l}
{\Xi_{1,0}^{\text{$\upgamma$}}(\text{t}),\Xi_{1,1}^{\text{$\upgamma$}}(\text{t}),\ldots,\Xi_{1,\text{M} - 1}^{\text{$\upgamma$}}(\text{t}),\Xi_{2,0}^{\text{$\upgamma$}}(\text{t}),\Xi_{2,1}^{\text{$\upgamma$}}(\text{t}),\ldots,\Xi_{2,\text{M} - 1}^{\text{$\upgamma$}}(\text{t}),\ldots,} \vspace{5pt}\\
{\Xi_{2^{\text{k} - 1},0}^{\text{$\upgamma$}}(\text{t}),\Xi_{2^{\text{k} - 1},1}^{\text{$\upgamma$}}(\text{t}),\ldots,\Xi_{2^{\text{k} - 1},\text{M}}^{\text{$\upgamma$}}(\text{t})} \\
\end{array} \right\rbrack^{\text{T}}
\tag{14}
\end{equation}

In the calculation process, we take $2^{\text{k} - 1}(\text{M} + 1) = \widetilde{\text{$\upsigma$}} $ which shows the total FOBW basis, and T represents the usual transpose.

The following section presents the FOBWs scheme for the evaluation of the variable-order fractional oscillation equation.


\section{Proposed Approach}
\label{sect:sec6-mathematics-2384216}

It has been revealed from the literature review that analysis of the physical, engineering, and mathematical models associated with variable-order fractional derivatives rather than derivatives of constant fractional order or integer order provides highly significant results. Therefore, motivated by the nice performance of the existing approaches given in Table~\ref{tabref:mathematics-2384216-t001}, we apply an effective wavelet approach for the numerical analysis and simulation of the variable-order nonlinear fractional oscillation equation.

The nonlinear model given in Equation (1) can be expressed as
\begin{equation}
\label{eq:FD21-mathematics-2384216}
\mathbb{Q}\left( {\text{D}_{0,\text{t}}^{\text{$\upalpha$}(\text{t})}\Im(\text{t}),\Im^{\prime}(\text{t}),\Im(\text{t}),\Phi\left( {\text{f},\text{$\upomega$},\text{t}} \right),\text{t}} \right) = 0
\tag{15}
\end{equation}
with the condition

\vspace{-9pt}
\begin{equation}
\label{eq:FD22-mathematics-2384216}
\Im(0) = 1,\quad\Im^{\prime}(0) = 0.
\tag{16}
\end{equation}

To determine the solutions of the above system, the procedure of the mentioned wavelet approach is provided stepwise as follows:

\textbf{\boldmath{Step I:}} The proposed method as well as the approximation through FOBWs totally depend on the range of $\text{$\upalpha$}(\text{t}). $ Since $\text{$\upalpha$}(\text{t}) \in \left( {1,2} \right\rbrack, $ then approximate the second-order derivative of an unknown function as a linear combination of truncated FOBWs using Equation (12) as

\vspace{-6pt}
\begin{equation}
\label{eq:FD23-mathematics-2384216}
\Im^{(2)}(\text{t}) \approx \text{U}_{}^{\text{T}}\text{$\Psi$}_{\widetilde{\text{$\upsigma$}}}^{\text{$\upgamma$}}(\text{t}),
\tag{17}
\end{equation}
where $\text{$\Psi$}_{\widetilde{\text{$\upsigma$}}}^{\text{$\upgamma$}}(\text{t}) $ is given in Equations (9) and (14) and $\text{U} $ is wavelet the coefficients vector.

Using Equations (3)--(5) on Equation (17), we get
\begin{equation}
\label{eq:FD24-mathematics-2384216}
\begin{array}{l}
{\text{I}_{0,\text{t}}^{2}\left( {\Im^{(2)}(\text{t})} \right) \approx \text{I}_{0,\text{t}}^{2}\left( {\text{U}_{}^{\text{T}}\text{$\Psi$}_{\widetilde{\text{$\upsigma$}}}^{\text{$\upgamma$}}(\text{t})} \right)}\vspace{5pt} \\
{\Im(\text{t}) - {\sum\limits_{\text{j} = 0}^{{\lceil 2\rceil} - 1}{\Im^{(2)}(0)\frac{\text{t}^{\text{j}}}{\text{j}!}}} = \text{U}_{}^{\text{T}}\left( {\text{I}_{0,\text{t}}^{2}\text{$\Psi$}_{\widetilde{\text{$\upsigma$}}}^{\text{$\upgamma$}}(\text{t})} \right)}\vspace{5pt} \\
{\Im(\text{t})\  - \text{t}\Im^{(1)}(0) - \Im(0) = \text{U}_{}^{\text{T}}\left( {\text{I}_{0,\text{t}}^{2}\text{$\Psi$}_{\widetilde{\text{$\upsigma$}}}^{\text{$\upgamma$}}(\text{t})} \right)} \vspace{5pt}\\
{\Im(\text{t})\  = \text{U}_{}^{\text{T}}\left( {\text{I}_{0,\text{t}}^{2}\text{$\Psi$}_{\widetilde{\text{$\upsigma$}}}^{\text{$\upgamma$}}(\text{t})} \right) + \Im(0) + \text{t}\Im^{(1)}(0),} \\
\end{array}
\tag{18}
\end{equation}
where $\text{I}_{0,\text{t}}^{2}\text{$\Psi$}_{\widetilde{\text{$\upsigma$}}}^{\text{$\upgamma$}}(\text{t}) $ is calculated directly by using Equation (3) on a known function $\text{$\Psi$}_{\widetilde{\text{$\upsigma$}}}^{\text{$\upgamma$}}(\text{t}) $ for different $\widetilde{\text{$\upsigma$}} $.

\textbf{\boldmath{Step II:}} Using Equation (5) with the range of $\text{$\upalpha$}(\text{t}) \in \left( {1,2} \right\rbrack, $ we get
\begin{equation}
\label{eq:FD25-mathematics-2384216}
\text{D}_{0,\text{t}}^{\text{$\upalpha$}(\text{t})}\Im(\text{t}) = \text{U}_{}^{\text{T}}\text{I}_{0,\text{t}}^{2 - \text{$\upalpha$}(\text{t})}\text{$\Psi$}_{\widetilde{\text{$\upsigma$}}}^{\text{$\upgamma$}}(\text{t}) + {\sum\limits_{\text{j} = {\lceil{\text{$\upalpha$}(\text{t})}\rceil}}^{1}{\Im^{(\text{j})}(0)\frac{\text{t}^{\text{j} - \text{$\upalpha$}(\text{t})}}{(\text{j} - \text{$\upalpha$}(\text{t}))!}}}.
\tag{19}
\end{equation}

\textbf{\boldmath{Step III:}} Substituting Equations (18) and (19) in the given system of Equation (15), we~get
\begin{equation}
\label{eq:FD26-mathematics-2384216}
\mathbb{Q}\left( \begin{array}{l}
{\text{U}_{}^{\text{T}}\text{I}_{0,\text{t}}^{2 - \text{$\upalpha$}(\text{t})}\text{$\Psi$}_{\widetilde{\text{$\upsigma$}}}^{\text{$\upgamma$}}(\text{t}) + {\sum\limits_{\text{j} = {\lceil{\text{$\upalpha$}(\text{t})}\rceil}}^{1}{\Im^{(\text{j})}(0)\frac{\text{t}^{\text{j} - \text{$\upalpha$}(\text{t})}}{(\text{j} - \text{$\upalpha$}(\text{t}))!}}},\ \text{U}_{}^{\text{T}}\text{I}_{0,\text{t}}^{}\text{$\Psi$}_{\widetilde{\text{$\upsigma$}}}^{\text{$\upgamma$}}(\text{t}) + \Im^{(1)}(0),} \vspace{5pt}\\
{\text{U}_{}^{\text{T}}\text{I}_{0,\text{t}}^{2}\text{$\Psi$}_{\widetilde{\text{$\upsigma$}}}^{\text{$\upgamma$}}(\text{t}) + \text{t}\Im^{(1)}(0) + \Im(0),\ \Phi\left( {\text{$\upomega$},\text{f},\text{t}} \right),\text{t}} \\
\end{array} \right) = 0.
\tag{20}
\end{equation}

\textbf{\boldmath{Step IV:}} The set of n non-linear algebraic equations is acquired via collocating the Equation (20) at appropriate Chebyshev grids $\text{t}_{\text{r}} $ as
\begin{equation}
\label{eq:FD27-mathematics-2384216}
\mathbb{Q}\left( \begin{array}{l}
{\text{U}_{}^{\text{T}}\text{I}_{0,\text{t}}^{2 - \text{$\upalpha$}(\text{t}_{\text{r}})}\text{$\Psi$}_{\widetilde{\text{$\upsigma$}}}^{\text{$\upgamma$}}(\text{t}_{\text{r}}) + {\sum\limits_{\text{j} = {\lceil{\text{$\upalpha$}(\text{t}_{\text{r}})}\rceil}}^{1}{\Im^{(\text{j})}(0)\frac{\text{t}_{\text{r}}{}^{\text{j} - \text{$\upalpha$}(\text{t}_{\text{r}})}}{(\text{j} - \text{$\upalpha$}(\text{t}_{\text{r}}))!}}},\ \text{U}_{}^{\text{T}}\text{I}_{0,\text{t}}^{}\text{$\Psi$}_{\widetilde{\text{$\upsigma$}}}^{\text{$\upgamma$}}(\text{t}_{\text{r}}) + \Im^{(1)}(0),} \vspace{5pt} \\
{\text{U}_{}^{\text{T}}\text{I}_{0,\text{t}}^{2}\text{$\Psi$}_{\widetilde{\text{$\upsigma$}}}^{\text{$\upgamma$}}(\text{t}_{\text{r}}) + \text{t}_{\text{r}}\Im^{(1)}(0) + \Im(0),\ \Phi\left( {\text{$\upomega$},\text{f},\text{t}_{\text{r}}} \right),\text{t}_{\text{r}}} \\
\end{array} \right) = 0,
\tag{21}
\end{equation}
where appropriate collocation grids $\text{t}_{\text{r}} $ is given by
\begin{equation}
\nonumber\label{eq:FD28-mathematics-2384216}
\text{t}_{\text{r}} = \frac{1}{2}\cos \left( \frac{(\text{r} - 1/2)\text{$\uppi$}}{\widetilde{\text{$\upsigma$}}} \right) + \frac{1}{2};\quad\text{r} = 1,2,\ldots,\widetilde{\text{$\upsigma$}}.
\end{equation}

\textbf{\boldmath{Step V:}} Solve the algebraic set of equations formed in Equation (21), we can easily find the unknown wavelet coefficient vectors $\text{U} $.

\textbf{\boldmath{Step VI:}} Using the value $\text{U} $ in Equation (18), we determine the wavelet approximation of $\Im(\text{t}). $

Procedure completed.

\begin{figure}[H]
\begin{adjustwidth}{-\extralength}{0cm}
\centering
\includegraphics[scale=1.1]{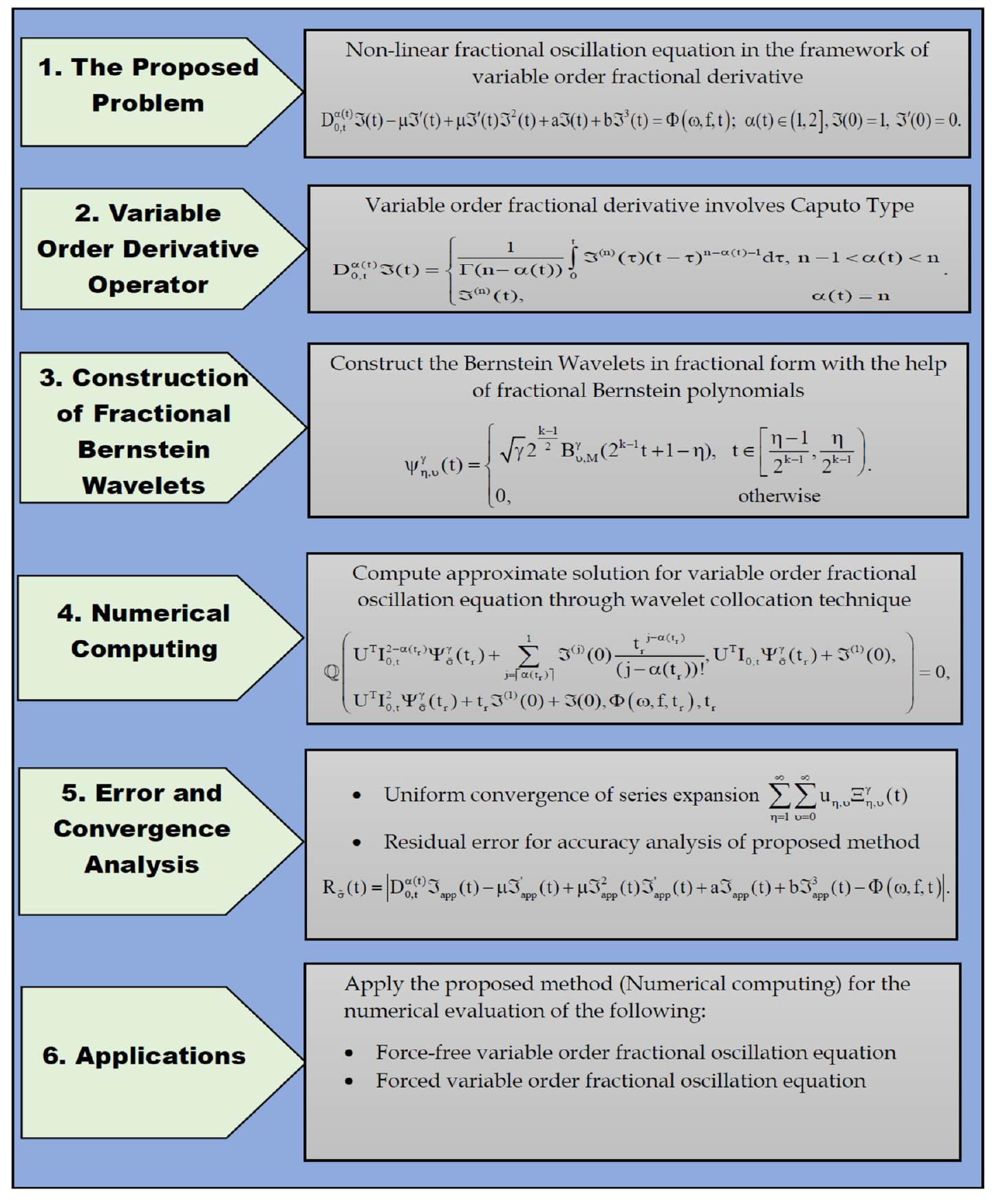}
\end{adjustwidth}
\caption{Graphical structure of Proposed Wavelets Technique.}
\label{fig:mathematics-2384216-f002}
\end{figure}

The graphical structure of the proposed scheme is represented in Figure~\ref{fig:mathematics-2384216-f002}.    


\section{Convergence Analysis}
\label{sect:sec7-mathematics-2384216}

The following theorem gives the uniform convergence of series expansion of function $\Im(\text{t}) \in \text{L}^{2}\lbrack 0,1) $.

\vspace{12pt}
\noindent \textbf{\boldmath{Theorem}} \textbf{\boldmath{1.}} \emph{Let} $\Im(\text{t}) \in \text{L}^{2}\lbrack 0,1) $ \emph{is a continuous bounded function and it can be approximated through the infinite series of FOBWs basis according to Equation (11), then this series converges uniformly to the function} $\Im(\text{t}) $. 

\vspace{12pt}
\noindent \textbf{\boldmath{Proof of Theorem}} \textbf{\boldmath{1.}} Consider the approximations of function $\Im(\text{t}) $ according to given definition of Equation (11) as
\begin{equation}
\label{eq:FD29-mathematics-2384216}
\Im(\text{t}) \approx {\sum\limits_{\text{$\upeta$} = 1}^{\infty}{\sum\limits_{\text{$\upupsilon$} = 0}^{\infty}{\text{u}_{\text{$\upeta$},\text{$\upupsilon$}}\Xi_{\text{$\upeta$},\text{$\upupsilon$}}^{\text{$\upgamma$}}(\text{t})}}},
\tag{22}
\end{equation}
then the wavelet coefficient is obtained by Equation (22) as
\begin{equation}
\label{eq:FD30-mathematics-2384216}
\begin{array}{cl}
\text{u}_{\text{$\upeta$},\text{$\upupsilon$}} & {= \left\langle {\Im(\text{t}),\Xi_{\text{$\upeta$},\text{$\upupsilon$}}^{\text{$\upgamma$}}(\text{t})} \right\rangle_{\text{$\Omega$}_{\text{k},\text{$\upeta$}}(\text{t})}} \vspace{6pt}\\
 & {= {\int\limits_{0}^{1}{\Im(\text{t})}}\Xi_{\text{$\upeta$},\text{$\upupsilon$}}^{\text{$\upgamma$}}(\text{t})\text{$\Omega$}_{\text{k},\text{$\upeta$}}(\text{t})\text{d}\text{t}.} \\
\end{array}
\tag{23}
\end{equation}

Using the definition of FOBWs in Equation (23), we obtain
\begin{equation}
\label{eq:FD31-mathematics-2384216}
\begin{array}{cl}
\text{u}_{\text{$\upeta$},\text{$\upupsilon$}} & {= {\int\limits_{\frac{\text{$\upeta$} - 1}{2^{\text{k} - 1}}}^{\frac{\text{$\upeta$}}{2^{\text{k} - 1}}}{\Im(\text{t})}}\sqrt{\text{$\upgamma$}}2^{(\text{k} - 1)/2}\text{B}_{\text{$\upupsilon$},\text{M}}^{\text{$\upgamma$}}(1 + 2^{\text{k} - 1}\text{t} - \text{$\upeta$})\text{d}\text{t}\ } \\
 & {= 2^{(\text{k} - 1)/2}\sqrt{\text{$\upgamma$}}\ {\int\limits_{\frac{\text{$\upeta$} - 1}{2^{\text{k} - 1}}}^{\frac{\text{$\upeta$}}{2^{\text{k} - 1}}}{\Im(\text{t})}}\text{B}_{\text{$\upupsilon$},\text{M}}^{\text{$\upgamma$}}(1 + 2^{\text{k} - 1}\text{t} - \text{$\upeta$})\text{d}\text{t}.\ } \\
\end{array}
\tag{24}
\end{equation}

Let $1 + 2^{\text{k} - 1}\text{t} - \text{$\upeta$} = \text{$\uptau$} $, then from Equation (24) we obtain
\begin{equation}
\label{eq:FD32-mathematics-2384216}
\begin{array}{cl}
\text{u}_{\text{$\upeta$},\text{$\upupsilon$}} & {= 2^{(\text{k} - 1)/2}\sqrt{\text{$\upgamma$}}\ {\int\limits_{0}^{1}{\Im\left( \frac{\text{$\uptau$} - 1 + \text{$\upeta$}}{2^{\text{k} - 1}} \right)}}\text{B}_{\text{$\upupsilon$},\text{M}}^{\text{$\upgamma$}}(\text{$\uptau$})\frac{\text{d}\text{$\uptau$}}{2^{\text{k} - 1}}} \vspace{6pt}\\
& {= \frac{1}{2^{(\text{k} - 1)/2}}\sqrt{\text{$\upgamma$}}\ {\int\limits_{0}^{1}{\Im\left( \frac{\text{$\uptau$} - 1 + \text{$\upeta$}}{2^{\text{k} - 1}} \right)}}\text{B}_{\text{$\upupsilon$},\text{M}}^{\text{$\upgamma$}}(\text{$\uptau$})\text{d}\text{$\uptau$},} \\
\end{array}
\tag{25}
\end{equation}
using the generalized mean value theorem of integrals in Equation (25), we get
\begin{equation}
\label{eq:FD33-mathematics-2384216}
\text{u}_{\text{$\upeta$},\text{$\upupsilon$}} = \frac{1}{2^{(\text{k} - 1)/2}}\sqrt{\text{$\upgamma$}}\ \Im\left( \frac{\text{z} - 1 + \text{$\upeta$}}{2^{\text{k} - 1}} \right){\int\limits_{0}^{1}{\text{B}_{\text{$\upupsilon$},\text{M}}^{\text{$\upgamma$}}(\text{$\uptau$})\text{d}\text{$\uptau$}}},\quad\text{z} \in \left( {0,1} \right).
\tag{26}
\end{equation}

Now by taking modulus both sides of Equation (26), we obtain
\begin{equation}
\label{eq:FD34-mathematics-2384216}
\left| \text{u}_{\text{$\upeta$},\text{$\upupsilon$}} \right| = \left| {\frac{1}{2^{(\text{k} - 1)/2}}\sqrt{\text{$\upgamma$}}} \right|\ \left| {\Im\left( \frac{\text{z} - 1 + \text{$\upeta$}}{2^{\text{k} - 1}} \right)} \right|\ {\int\limits_{0}^{1}{\left| {\text{B}_{\text{$\upupsilon$},\text{M}}^{\text{$\upgamma$}}(\text{$\uptau$})} \right|\text{d}\text{$\uptau$}}}.
\tag{27}
\end{equation}

Since $\Im(\text{t}) $ is a continuous bounded function, then $\left| {\Im\left( \frac{\text{z} - 1 + \text{$\upeta$}}{2^{\text{k} - 1}} \right)} \right| \leq \text{${\uprho}$}. $

Then Equation (27) implies
\begin{equation}
\label{eq:FD35-mathematics-2384216}
\left| \text{u}_{\text{$\upeta$},\text{$\upupsilon$}} \right|\  \leq \text{${\uprho}$}\frac{\sqrt{\text{$\upgamma$}}}{2^{(\text{k} - 1)/2}}{\int\limits_{0}^{1}{\left| {\text{B}_{\text{$\upupsilon$},\text{M}}^{\text{$\upgamma$}}(\text{$\uptau$})} \right|\text{d}\text{$\uptau$}}}.
\tag{28}
\end{equation}

The fractional order Bernstein polynomials given in Equation (8) are given as
\begin{equation}
\label{eq:FD36-mathematics-2384216}
\text{B}_{\text{$\upupsilon$},\text{M}}^{\text{$\upgamma$}}(\text{t}) = \sqrt{1 + 2\text{M} - 2\text{$\upupsilon$}}\quad{\sum\limits_{\text{i} = 0}^{\text{$\upupsilon$}}\left( {- 1} \right)^{\text{i}}}\frac{\begin{pmatrix}
{1 + 2\text{M} - \text{i}} \\
{\text{$\upupsilon$} - \text{i}} \\
\end{pmatrix}\begin{pmatrix}
\text{$\upupsilon$} \\
\text{i} \\
\end{pmatrix}}{\begin{pmatrix}
{\text{M} - \text{i}} \\
{\text{$\upupsilon$} - \text{i}} \\
\end{pmatrix}}{\widetilde{\text{B}}}_{\text{$\upupsilon$} - \text{i},\text{M} - \text{i}}^{\text{$\upgamma$}}(\text{t}).
\tag{29}
\end{equation}

Therefore, Equation (29) implies
  
\begin{adjustwidth}{-\extralength}{0cm}
\centering 
\begin{equation}
\nonumber\label{eq:FD37-mathematics-2384216}
\begin{array}{cl}
{\left| {\text{B}_{\text{$\upupsilon$},\text{M}}^{\text{$\upgamma$}}(\text{t})} \right| } & {= \sqrt{1 + 2\text{M} - 2\text{$\upupsilon$}}\ {\sum\limits_{\text{i} = 0}^{\text{$\upupsilon$}}\left( {- 1} \right)^{\text{i}}}\frac{\left( \begin{array}{c}
{1 + 2\text{M} - \text{i}} \\
{\text{$\upupsilon$} - \text{i}} \\
\end{array} \right)\left( \begin{array}{c}
\text{$\upupsilon$} \\
\text{i} \\
\end{array} \right)}{\left( \begin{array}{c}
{\text{M} - \text{i}} \\
{\text{$\upupsilon$} - \text{i}} \\
\end{array} \right)}\left| {{\widetilde{\text{B}}}_{\text{$\upupsilon$} - \text{i},\text{M} - \text{i}}^{\text{$\upgamma$}}(\text{t})} \right|} \vspace{6pt}\\
 & {\leq \sqrt{1 + 2\text{M} - 2\text{$\upupsilon$}}\ {\sum\limits_{\text{i} = 0}^{\text{$\upupsilon$}}\frac{\left( \begin{array}{c}
{1 + 2\text{M} - \text{i}} \\
{\text{$\upupsilon$} - \text{i}} \\
\end{array} \right)\left( \begin{array}{c}
\text{$\upupsilon$} \\
\text{i} \\
\end{array} \right)}{\left( \begin{array}{c}
{\text{M} - \text{i}} \\
{\text{$\upupsilon$} - \text{i}} \\
\end{array} \right)}} \times {\sum\limits_{\text{r} = 0}^{\text{M} - \text{$\upupsilon$}}{\left( \begin{array}{c}
{\text{M} - \text{i}} \\
{\text{$\upupsilon$} - \text{i}} \\
\end{array} \right)\left( \begin{array}{c}
{\text{M} - \text{$\upupsilon$}} \\
\text{r} \\
\end{array} \right)\left| \text{t}^{\text{$\upgamma$}(\text{$\upupsilon$} - \text{i} + \text{r})} \right|}}.} \\
\end{array}
\end{equation}
\end{adjustwidth}

Using the above result in Equation (28), we get
\begin{equation}
\label{eq:FD38-mathematics-2384216}
\begin{array}{ll}
\left| \text{u}_{\text{$\upeta$},\text{$\upupsilon$}} \right| & {\leq \text{${\uprho}$}\frac{\sqrt{\text{$\upgamma$}}}{2^{(\text{k} - 1)/2}}{\int\limits_{0}^{1}{\sqrt{1 + 2\text{M} - 2\text{$\upupsilon$}}\ {\sum\limits_{\text{i} = 0}^{\text{$\upupsilon$}}\frac{\left( \begin{array}{c}
{1 + 2\text{M} - \text{i}} \\
{\text{$\upupsilon$} - \text{i}} \\
\end{array} \right)\left( \begin{array}{c}
\text{$\upupsilon$} \\
\text{i} \\
\end{array} \right)}{\left( \begin{array}{c}
{\text{M} - \text{i}} \\
{\text{$\upupsilon$} - \text{i}} \\
\end{array} \right)}}}}} \vspace{4pt}\\
 & {\times {\sum\limits_{\text{r} = 0}^{\text{M} - \text{$\upupsilon$}}{\left( \begin{array}{c}
{\text{M} - \text{i}} \\
{\text{$\upupsilon$} - \text{i}} \\
\end{array} \right)\left( \begin{array}{c}
{\text{M} - \text{$\upupsilon$}} \\
\text{r} \\
\end{array} \right)\left| \text{t}^{\text{$\upgamma$}(\text{$\upupsilon$} - \text{i} + \text{r})} \right|}}\text{d}\text{$\uptau$}} \vspace{8pt}\\
 & {\leq \text{${\uprho}$}\frac{\sqrt{\text{$\upgamma$}}}{2^{(\text{k} - 1)/2}}\sqrt{1 + 2\text{M} - 2\text{$\upupsilon$}}{\sum\limits_{\text{i} = 0}^{\text{$\upupsilon$}}{\left( \begin{array}{c}
{1 + 2\text{M} - \text{i}} \\
{\text{$\upupsilon$} - \text{i}} \\
\end{array} \right)\left( \begin{array}{c}
\text{$\upupsilon$} \\
\text{i} \\
\end{array} \right)}} \times {\sum\limits_{\text{r} = 0}^{\text{M} - \text{$\upupsilon$}}\left( \begin{array}{c}
{\text{M} - \text{$\upupsilon$}} \\
\text{r} \\
\end{array} \right)}} \vspace{8pt}\\
 & {\leq \text{${\uprho}$}\frac{\sqrt{\text{$\upgamma$}}}{2^{(\text{k} - 1)/2}}2^{\text{M} - \text{$\upupsilon$}}\sqrt{1 + 2\text{M} - 2\text{$\upupsilon$}}\left( \begin{array}{c}
{1 + 2\text{M} - \text{i} + \text{$\upupsilon$}} \\
\text{$\upupsilon$} \\
\end{array} \right).} \\
\end{array}
\tag{30}
\end{equation}

Since $\Im(\text{t}) $ is bounded and continuous and $\left| \text{u}_{\text{$\upeta$},\text{$\upupsilon$}} \right|\  $ is finite for the existing parameters, therefore $\sum\limits_{\text{$\upeta$} = 1}^{\infty}{\sum\limits_{\text{$\upupsilon$} = 0}^{\infty}\text{u}_{\text{$\upeta$},\text{$\upupsilon$}}} $ is absolutely convergent by definition of convergence of series.

Hence the series expansion $\sum\limits_{\text{$\upeta$} = 1}^{\infty}{\sum\limits_{\text{$\upupsilon$} = 0}^{\infty}{\text{u}_{\text{$\upeta$},\text{$\upupsilon$}}\Xi_{\text{$\upeta$},\text{$\upupsilon$}}^{\text{$\upgamma$}}(\text{t})}} $ convergence uniformly to $\Im(\text{t}) $. $\square$

\vspace{12pt}
\noindent \textbf{\boldmath{Theorem}} \textls[-25]{\textbf{\boldmath{2.}} \emph{Let} $\Im(\text{t}) \in \text{L}^{2}\lbrack 0,1) $ \emph{is a continuous bounded function and} $\Im_{\widetilde{\text{$\upsigma$}}}(\text{t}) = {\sum\limits_{\text{$\upeta$} = 1}^{2^{\text{k} - 1}}{\sum\limits_{\text{$\upupsilon$} = 0}^{\text{M}}{\text{u}_{\text{$\upeta$},\text{$\upupsilon$}}\Xi_{\text{$\upeta$},\text{$\upupsilon$}}^{\text{$\upgamma$}}(\text{t})}}} $} \emph{be a FOBWs approximation of $\Im(\text{t}), $ then the upper bound of error is estimated as}\begin{equation}
\nonumber\label{eq:FD39-mathematics-2384216}
\left\| {\Im(\text{t}) - \Im_{\widetilde{\text{$\upsigma$}}}(\text{t})} \right\|_{\text{L}_{\text{$\upomega$}}^{2}\lbrack 0,1\rbrack} \leq \text{$\upkappa$},
\end{equation}
\emph{where}

\vspace{-3pt}
\begin{equation}
\label{eq:FD40-mathematics-2384216}
\text{$\upkappa$} = \left( {\sum\limits_{\text{$\upeta$} = 2^{\text{k} - 1} + 1}^{\infty}{\sum\limits_{\text{$\upupsilon$} = \text{M}}^{\infty}\left| \text{u}_{\text{$\upeta$},\text{$\upupsilon$}} \right|^{2}}} \right)^{1/2},
\tag{31}
\end{equation}
 \emph{and} $\text{u}_{\text{$\upeta$},\text{$\upupsilon$}} $ \emph{is given in Equation (13).}

\vspace{12pt}
\noindent \textbf{\boldmath{Proof of Theorem}} \textbf{\boldmath{2.}} Since $\Im_{\widetilde{\text{$\upsigma$}}}(\text{t}) $ be the FOBW’s approximation of $\Im(\text{t}), $ then

\begin{adjustwidth}{-\extralength}{0cm}
\centering 
\begin{equation}
\label{eq:FD41-mathematics-2384216}
\begin{array}{cl}
\left\| {\Im(\text{t}) - \Im_{\widetilde{\text{$\upsigma$}}}(\text{t})} \right\|_{\text{L}_{\text{$\upomega$}}^{2}\lbrack 0,1\rbrack}^{2} & {= \left\| {\Im(\text{t}) - {\sum\limits_{\text{$\upeta$} = 1}^{2^{\text{k} - 1}}{\sum\limits_{\text{$\upupsilon$} = 0}^{\text{M}}{\text{u}_{\text{$\upeta$},\text{$\upupsilon$}}\Xi_{\text{$\upeta$},\text{$\upupsilon$}}^{\text{$\upgamma$}}(\text{t})}}}}
\right\|_{\text{L}_{\text{$\Omega$}}^{2}\lbrack 0,1\rbrack}^{2}} \vspace{6pt}\\
& {= {\int\limits_{0}^{1}{\left| {\Im(\text{t}) - {\sum\limits_{\text{$\upeta$} = 1}^{2^{\text{k} - 1}}{\sum\limits_{\text{$\upupsilon$} = 0}^{\text{M}}{\text{u}_{\text{$\upeta$},\text{$\upupsilon$}}\Xi_{\text{$\upeta$},\text{$\upupsilon$}}^{\text{$\upgamma$}}(\text{t})}}}} \right|^{2}\ \text{$\Omega$}_{\text{k},\text{$\upeta$}}(\text{t})\text{d}\text{t}}}\ } \vspace{6pt}\\
& {= {\int\limits_{0}^{1}{\left| {{\sum\limits_{\text{$\upeta$} = 1}^{\infty}{\sum\limits_{\text{$\upupsilon$} = 0}^{\infty}{\text{u}_{\text{$\upeta$},\text{$\upupsilon$}}\Xi_{\text{$\upeta$},\text{$\upupsilon$}}^{\text{$\upgamma$}}(\text{t})}}} - {\sum\limits_{\text{$\upeta$} = 1}^{2^{\text{k} - 1}}{\sum\limits_{\text{$\upupsilon$} = 0}^{\text{M}}{\text{u}_{\text{$\upeta$},\text{$\upupsilon$}}\Xi_{\text{$\upeta$},\text{$\upupsilon$}}^{\text{$\upgamma$}}(\text{t})}}}} \right|^{2}\ \text{$\Omega$}_{\text{k},\text{$\upeta$}}(\text{t})\text{d}\text{t}}}} \vspace{6pt}\\
& {= {\int\limits_{0}^{1}{\left| {\sum\limits_{\text{$\upeta$} = 2^{\text{k} - 1} + 1}^{\infty}{\sum\limits_{\text{$\upupsilon$} = \text{M} + 1}^{\infty}{\text{u}_{\text{$\upeta$},\text{$\upupsilon$}}\Xi_{\text{$\upeta$},\text{$\upupsilon$}}^{\text{$\upgamma$}}(\text{t})}}} \right|^{2}\ \text{$\Omega$}_{\text{k},\text{$\upeta$}}(\text{t})\text{d}\text{t}}}} \vspace{6pt}\\
& {= {\sum\limits_{\text{$\upeta$} = 2^{\text{k} - 1} + 1}^{\infty}{\sum\limits_{\text{$\upupsilon$} = \text{M}}^{\infty}\left| \text{u}_{\text{$\upeta$},\text{$\upupsilon$}}
\right|^{2}}}{\int\limits_{0}^{1}{\Xi_{\text{$\upeta$},\text{$\upupsilon$}}^{\text{$\upgamma$}}(\text{t})\Xi_{\text{$\upeta$},\text{$\upupsilon$}}^{\text{$\upgamma$}}(\text{t})
\ \text{$\Omega$}_{\text{k},\text{$\upeta$}}(\text{t})\text{d}\text{t}}}.} \\
\end{array}
\tag{32}
\end{equation}
\end{adjustwidth}

Due to the orthonormality of FOBWs, from Equation (32) we obtain:\begin{equation}
\label{eq:FD42-mathematics-2384216}
\left\| {\Im(\text{t}) - \Im_{\widetilde{\text{$\upsigma$}}}(\text{t})} \right\|_{\text{L}_{\text{$\upomega$}}^{2}\lbrack 0,1\rbrack}^{2} = {\sum\limits_{\text{$\upeta$} = 2^{\text{k} - 1} + 1}^{\infty}{\sum\limits_{\text{$\upupsilon$} = \text{M} + 1}^{\infty}\left| \text{u}_{\text{$\upeta$},\text{$\upupsilon$}} \right|^{2}}},
\tag{33}
\end{equation}
where $\text{u}_{\text{$\upeta$},\text{$\upupsilon$}} $ is given in Equation (13).

By taking square-roots, we get
\begin{equation}
\label{eq:FD43-mathematics-2384216}
\left\| {\Im(\text{t}) - \Im_{\widetilde{\text{$\upsigma$}}}(\text{t})} \right\|_{\text{L}_{\text{$\upomega$}}^{2}\lbrack 0,1\rbrack}^{} = \left( {\sum\limits_{\text{$\upeta$} = 2^{\text{k} - 1} + 1}^{\infty}{\sum\limits_{\text{$\upupsilon$} = \text{M} + 1}^{\infty}\left| \text{u}_{\text{$\upeta$},\text{$\upupsilon$}} \right|^{2}}} \right)^{1/2}.
\tag{34}
\end{equation}

And from Equation (30) of Theorem 1, we have
\begin{equation}
\label{eq:FD44-mathematics-2384216}
\left| \text{u}_{\text{$\upeta$},\text{$\upupsilon$}} \right|\  \leq \text{${\uprho}$}\sqrt{\text{$\upgamma$}}\frac{2^{\text{M} - \text{$\upupsilon$}}}{2^{(\text{k} - 1)/2}}\sqrt{1 + 2\text{M} - 2\text{$\upupsilon$}}\ \begin{pmatrix}
{1 + 2\text{M} - \text{i} + \text{$\upupsilon$}} \\
\text{$\upupsilon$} \\
\end{pmatrix}.
\tag{35}
\end{equation}

Therefore, from Equations (33)--(35), we have
\begin{equation}
\nonumber\label{eq:FD45-mathematics-2384216}
\left\| {\Im(\text{t}) - \Im_{\widetilde{\text{$\upsigma$}}}(\text{t})} \right\|_{\text{L}_{\text{$\upomega$}}^{2}\lbrack 0,1\rbrack}^{} \leq \text{$\upkappa$},
\end{equation}
where $\text{$\upkappa$} $ is given in Equation (31).

Hence, the proof is complete. $\square$

\vspace{12pt}


\section{Numerical Examples}
\label{sect:sec8-mathematics-2384216}

The suggested approach is applied to the mentioned model (force-free and forced oscillation equations) to examine the performance of the approach for different parameters. All calculations are computed by the software Mathematica 7. In the examples, solutions are computed from t = 0 to t = 1 according to the parameters considered in the FOBWs basis. We can consider values t \textgreater{} 1 by modifying the range of FOBWs. The formula for absolute error is given for comparison purposes and to examine the efficiency of the mentioned~approach.
\begin{enumerate}[label=$\bullet$]
\item The absolute errors (AEs) between the wavelet approximation function $\Im_{\text{a}\text{p}\text{p}}^{}(\text{t}) $ and the analytical function $\Im(\text{t}) $ is computed as
\begin{equation}
\nonumber\label{eq:FD46-mathematics-2384216}
\left. \text{E}_{\text{A}\text{b}\text{s}}\ (\text{t}) = \  \middle| \Im(\text{t}) - \Im_{\text{a}\text{p}\text{p}}^{}(\text{t}) \middle| , \right.
\end{equation}
and the maximum absolute error (MAE) in this case is calculated as \linebreak $\text{M}\text{A}\text{E}\ \{\Im,\Im_{\text{a}\text{p}\text{p}}\} = \ \underset{\text{t} = \lbrack 0,1\rbrack}{\text{max}}|\Im(\text{t}) - \Im_{\text{a}\text{p}\text{p}}(\text{t})|. $
\item Since the analytical solutions of this model for fractional random order are not available, a residual error function $\text{R}_{\widetilde{\text{$\upsigma$}}}^{}(\text{t}) $ is introduced to measure the accuracy of the proposed approach as follows:

\begin{adjustwidth}{-\extralength}{0cm}
\begin{equation}
\nonumber\label{eq:FD47-mathematics-2384216}
\text{R}_{\widetilde{\text{$\upsigma$}}}^{}(\text{t}) = \left| {\text{D}_{0,\text{t}}^{\text{$\upalpha$}(\text{t})}\Im_{\text{a}\text{p}\text{p}}^{}(\text{t}) - \text{$\upmu$}\Im_{\text{a}\text{p}\text{p}}^{’}(\text{t}) + \text{$\upmu$}\Im_{\text{a}\text{p}\text{p}}^{2}(\text{t})\Im_{\text{a}\text{p}\text{p}}^{’}(\text{t}) + \text{a}\Im_{\text{a}\text{p}\text{p}}^{}(\text{t}) + \text{b}\Im_{\text{a}\text{p}\text{p}}^{3}(\text{t}) - \Phi\left( {\text{$\upomega$},\text{f},\text{t}} \right)} \right|.
\end{equation}
\end{adjustwidth}
\end{enumerate}

\vspace{12pt}
\noindent \textbf{\boldmath{Example}} \textbf{\boldmath{1.}} \emph{Consider the variable fractional order forced Duffing-Vander pol oscillator equation by replacing the constant fractional order~\cite{B29-mathematics-2384216} as}
\begin{equation}
\label{eq:FD48-mathematics-2384216}
\text{D}_{0,\text{t}}^{\text{$\upalpha$}(\text{t})}\Im(\text{t}) - \text{$\upmu$}\Im^{\prime}(\text{t}) + \text{$\upmu$}\Im^{\prime}(\text{t})\Im^{2}(\text{t}) + \text{a}\Im(\text{t}) + \text{b}\Im^{3}(\text{t}) = \text{f}\cos (\text{$\upomega$}\text{t}),\ \text{$\upalpha$}(\text{t}) \in \left( {1,2} \right\rbrack
\tag{36}
\end{equation}
\emph{with the initial value conditions}
\begin{equation}
\nonumber\label{eq:FD49-mathematics-2384216}
\Im(0) = 1,\ \Im^{\prime}(0) = 0.
\end{equation}

\emph{We solve the example for} $\widetilde{\text{$\upsigma$}} = 4,6\ (\text{k} = 1,\ \text{M} = 3,5)$ 
\emph{by mentioned scheme and simulate the model for different physically fascinating situations 
(single-well, double-well, and double-hump well) of the forced Duffing--Vander pol oscillator equation.}

\vspace{12pt}
In considering the problem, the following two cases of fractional order are considered:
\begin{enumerate}
\item[(i)] Constant order: $\text{$\upalpha$}(\text{t}) = 1.2,1.4,1.5,1.6,1.8$;
\item[(ii)] Variable order: $\text{$\upalpha$}(\text{t}) = 1 + \sin \text{t}$.
\end{enumerate}

\noindent\textbf{\boldmath{Physically fascinating conditions:}}

\begin{enumerate}
\item[(A)] \textbf{\boldmath{For Single-well}} $(\text{a},\text{b} > 0)$.
\end{enumerate}

\begin{figure}[H]
\includegraphics[scale=.88]{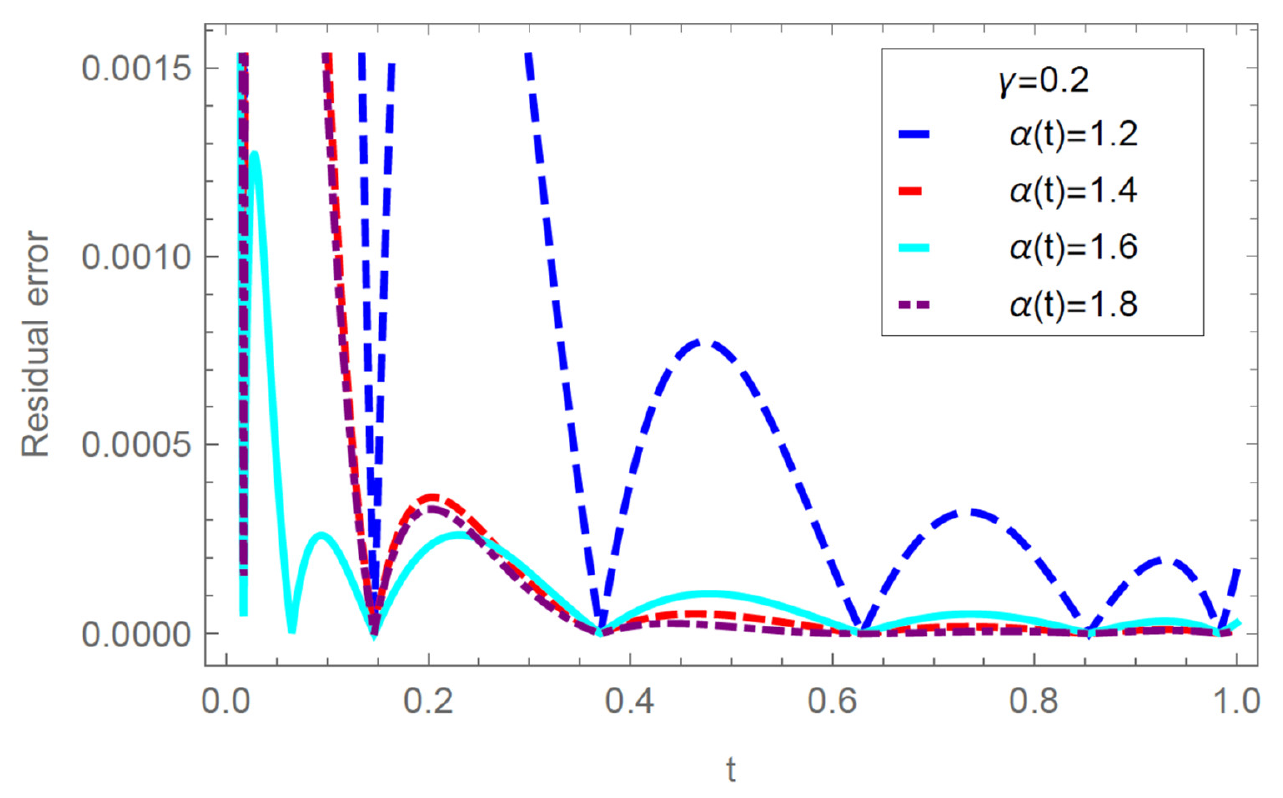}
\caption{Graph of Residual errors for single well case with k = 1, M = 5.}
\label{fig:mathematics-2384216-f003}
\end{figure}

The estimated AEs in the solutions of $\Im(\text{t}) $ with the comparison of the Legendre Wavelet--Picard scheme (LWPS) and the ultraspherical wavelets scheme (UWS) for $\text{$\upalpha$}(\text{t}) = 2 $ and different FOBWs bases are listed in Table~\ref{tabref:mathematics-2384216-t002}. It can be easily analyzed from Table~\ref{tabref:mathematics-2384216-t002} that the suitable value of $\text{$\upgamma$} $ is 1 for achieving the best accuracy in the solution of the given model, with $\text{$\upalpha$}(\text{t}) = 2 $ and the proposed approach is superior to UWS~\cite{B28-mathematics-2384216} and LWPS~\cite{B29-mathematics-2384216} by considering the RK-4 solution~\cite{B28-mathematics-2384216} as an approximated analytical solution. The residual errors in $\Im(\text{t}) $ for $\text{$\upalpha$}(\text{t}) = 1.5 $ and $\text{$\upalpha$}(\text{t}) = 1 + \sin \text{t} $ under different parameters mentioned are presented 
in Tables~\ref{tabref:mathematics-2384216-t003} and \ref{tabref:mathematics-2384216-t004}, respectively. In addition, the estimated residual errors in the solutions for $\text{$\upgamma$} = 0.2 $ and different selections of $\text{$\upalpha$}(\text{t}) $ are given in Table~\ref{tabref:mathematics-2384216-t005}. The graphical interpretation of residual errors of solutions for the single-well case with selected values of $\text{$\upalpha$}(\text{t}) $, and $\text{$\upgamma$} = 0.2 $ is shown in Figure~\ref{fig:mathematics-2384216-f003}. The computed solutions are obtained for the first time with the variable order of the introduced model in terms of residual errors.    
\begin{table}[H]
\tablesize{\small}
\caption{Estimated AEs for $\text{$\upalpha$}(\text{t}) = 2 $ and selected $\text{$\upgamma$} $ in Example 1.}
\label{tabref:mathematics-2384216-t002}
    
\setlength{\cellWidtha}{\textwidth/6-2\tabcolsep-0in}
\setlength{\cellWidthb}{\textwidth/6-2\tabcolsep-0in}
\setlength{\cellWidthc}{\textwidth/6-2\tabcolsep-0in}
\setlength{\cellWidthd}{\textwidth/6-2\tabcolsep-0in}
\setlength{\cellWidthe}{\textwidth/6-2\tabcolsep-0in}
\setlength{\cellWidthf}{\textwidth/6-2\tabcolsep-0in}
\scalebox{1}[1]{
\begin{tabularx}{\textwidth}{>{\centering\arraybackslash}m{\cellWidtha}>{\centering\arraybackslash}m{\cellWidthb}>{\centering\arraybackslash}m{\cellWidthc}>{\centering\arraybackslash}m{\cellWidthd}>{\centering\arraybackslash}m{\cellWidthe}>{\centering\arraybackslash}m{\cellWidthf}}
\toprule
\multicolumn{6}{>{\centering\arraybackslash}m{\cellWidtha + \cellWidthb + \cellWidthc + \cellWidthd + \cellWidthe + \cellWidthf+10\tabcolsep}}{\textbf{$\textbf{a}\ \textbf{=}\ \textbf{b}\ \textbf{=}\ \textbf{0.5}\textbf{,}\ \textbf{f}\ \textbf{=}\ \textbf{0.5}\textbf{,}\ \bm{\upmu}\ \textbf{=}\ \textbf{0.1}\textbf{,}\ \bm{\upomega}\ \textbf{=}\ \textbf{0.79} $}}\\
\cmidrule{1-6}
\multirow{2}{*}{\parbox{\cellWidtha}{\centering \textbf{t}}\vspace{-4pt}} & \multicolumn{3}{>{\centering\arraybackslash}m{\cellWidthb + \cellWidthc + \cellWidthd+4\tabcolsep}}{\textbf{
Proposed Approach, $\textbf{k}\ \textbf{=}\ \textbf{1}\textbf{,}\ \textbf{M}\ \textbf{=}\ \textbf{5} $}} & \multicolumn{2}{>{\centering\arraybackslash}m{\cellWidthe + \cellWidthf
+2\tabcolsep}}{\textbf{Reference Approach, $\textbf{M}\ \textbf{=}\ \textbf{6} $ }}\\
\cmidrule{2-6}
& \textbf{$\bm{\upgamma}\ \textbf{=}\ \textbf{0.5} $} & \textbf{$\bm{\upgamma}\ \textbf{=}\ \textbf{0.9} $} 
& \textbf{$\bm{\upgamma}\ \textbf{=}\ \textbf{1} $} & \textbf{UWS {\cite{B28-mathematics-2384216}} } 
& \textbf{LWPS {\cite{B29-mathematics-2384216}} }\\
\cmidrule{1-6}
0.1 & $4.2 \times 10^{- 7} $ & $1.1 \times 10^{- 7} $ & $1.2 \times 10^{- 7} $ & $4.0 \times 10^{- 8} $ & $2.1 \times 10^{- 8} $\\
\cmidrule{1-6}
0.3 & $1.7 \times 10^{- 6} $ & $5.3 \times 10^{- 8} $ & $1.1 \times 10^{- 7} $ & $1.2 \times 10^{- 7} $ & $5.2 \times 10^{- 8} $\\
\cmidrule{1-6}
0.5 & $1.8 \times 10^{- 6} $ & $1.5 \times 10^{- 7} $ & $1.3 \times 10^{- 7} $ & $6.1 \times 10^{- 7} $ & $2.8 \times 10^{- 7} $\\
\cmidrule{1-6}
0.7 & $2.6 \times 10^{- 6} $ & $2.9 \times 10^{- 7} $ & $3.1 \times 10^{- 8} $ & $1.6 \times 10^{- 6} $ & $1.4 \times 10^{- 6} $\\
\cmidrule{1-6}
0.9 & $2.9 \times 10^{- 6} $ & $3.7 \times 10^{- 7} $ & $2.8 \times 10^{- 8} $ & $3.3 \times 10^{- 6} $ & $2.5 \times 10^{- 6} $\\
\bottomrule
\end{tabularx}}
\end{table}
\vspace{-12pt}
    
\begin{table}[H]
\tablesize{\small}
\caption{\textls[-15]{Estimated residual errors of solutions for $\text{$\upalpha$}(\text{t}) = 1.5, $ $\text{k} = 1,\text{M} = 5, $ and selected $\text{$\upgamma$} $ in Example 1.}}
\label{tabref:mathematics-2384216-t003}
\setlength{\cellWidtha}{\textwidth/7-2\tabcolsep-0in}
\setlength{\cellWidthb}{\textwidth/7-2\tabcolsep-0in}
\setlength{\cellWidthc}{\textwidth/7-2\tabcolsep-0in}
\setlength{\cellWidthd}{\textwidth/7-2\tabcolsep-0in}
\setlength{\cellWidthe}{\textwidth/7-2\tabcolsep-0in}
\setlength{\cellWidthf}{\textwidth/7-2\tabcolsep-0in}
\setlength{\cellWidthg}{\textwidth/7-2\tabcolsep-0in}
\scalebox{1}[1]{
\begin{tabularx}{\textwidth}{>{\centering\arraybackslash}m{\cellWidtha}>{\centering\arraybackslash}m{\cellWidthb}>{\centering\arraybackslash}m{\cellWidthc}>{\centering\arraybackslash}m{\cellWidthd}>{\centering\arraybackslash}m{\cellWidthe}>{\centering\arraybackslash}m{\cellWidthf}>{\centering\arraybackslash}m{\cellWidthg}}
\toprule
\multicolumn{7}{>{\centering\arraybackslash}m{\cellWidtha + \cellWidthb + \cellWidthc + \cellWidthd + \cellWidthe + \cellWidthf + \cellWidthg+12\tabcolsep}}{\textbf{$\textbf{a}\ \textbf{=}\ \textbf{b}\ \textbf{=}\ \textbf{0.5}\textbf{,}\ \textbf{f}\ \textbf{=}\ \textbf{0.5}\textbf{,}\ \bm{\upmu}\ \textbf{=}\ \textbf{0.1}\textbf{,}\ \bm{\upomega}\ \textbf{=}\ \textbf{0.79} $}}\\
\cmidrule{1-7}
\textbf{t} & \textbf{$\bm{\upgamma}\ \textbf{=}\ \textbf{0.1} $} & \textbf{$\bm{\upgamma}\ \textbf{=}\ \textbf{0.2} $} & \textbf{$\bm{\upgamma}\ \textbf{=}\ \textbf{0.3} $} & \textbf{$\bm{\upgamma}\ \textbf{=}\ \textbf{0.5} $} & \textbf{$\bm{\upgamma}\ \textbf{=}\ \textbf{0.9} $} & \textbf{$\bm{\upgamma}\ \textbf{=}\ \textbf{1.0} $}\\
\cmidrule{1-7}
0.1 & $1.7 \times 10^{- 3} $ & $1.1 \times 10^{- 4} $ & $3.5 \times 10^{- 3} $ & $1.5 \times 10^{- 2} $ & $6.8 \times 10^{- 2} $ & $8.7 \times 10^{- 2} $\\
\cmidrule{1-7}
0.3 & $3.7 \times 10^{- 4} $ & $1.1 \times 10^{- 4} $ & $4.8 \times 10^{- 4} $ & $3.8 \times 10^{- 3} $ & $3.6 \times 10^{- 2} $ & $5.6 \times 10^{- 2} $\\
\cmidrule{1-7}
0.5 & $1.8 \times 10^{- 4} $ & $7.9 \times 10^{- 5} $ & $2.0 \times 10^{- 4} $ & $2.1 \times 10^{- 3} $ & $3.3 \times 10^{- 2} $ & $5.6 \times 10^{- 2} $\\
\cmidrule{1-7}
0.7 & $7.0 \times 10^{- 5} $ & $3.4 \times 10^{- 5} $ & $7.3 \times 10^{- 5} $ & $9.6 \times 10^{- 4} $ & $2.0 \times 10^{- 2} $ & $3.7 \times 10^{- 2} $\\
\cmidrule{1-7}
0.9 & $3.8 \times 10^{- 5} $ & $1.9 \times 10^{- 5} $ & $3.9 \times 10^{- 5} $ & $5.8 \times 10^{- 4} $ & $1.6 \times 10^{- 2} $ & $3.2 \times 10^{- 2} $\\
\bottomrule
\end{tabularx}}
\end{table}
\vspace{-12pt}
\begin{table}[H]
\tablesize{\small}
\caption{Estimated residual errors of solutions for $\text{$\upalpha$}(\text{t}) = 1 + \sin \text{t}, $ $\text{k} = 1,\text{M} = 5, $ 
and selected $\text{$\upgamma$} $ in Example 1.}
\label{tabref:mathematics-2384216-t004}
\setlength{\cellWidtha}{\textwidth/7-2\tabcolsep-0in}
\setlength{\cellWidthb}{\textwidth/7-2\tabcolsep-0in}
\setlength{\cellWidthc}{\textwidth/7-2\tabcolsep-0in}
\setlength{\cellWidthd}{\textwidth/7-2\tabcolsep-0in}
\setlength{\cellWidthe}{\textwidth/7-2\tabcolsep-0in}
\setlength{\cellWidthf}{\textwidth/7-2\tabcolsep-0in}
\setlength{\cellWidthg}{\textwidth/7-2\tabcolsep-0in}
\scalebox{1}[1]{
\begin{tabularx}{\textwidth}{>{\centering\arraybackslash}m{\cellWidtha}>{\centering\arraybackslash}m{\cellWidthb}>{\centering\arraybackslash}m{\cellWidthc}>{\centering\arraybackslash}m{\cellWidthd}>{\centering\arraybackslash}m{\cellWidthe}>{\centering\arraybackslash}m{\cellWidthf}>{\centering\arraybackslash}m{\cellWidthg}}
\toprule
\multicolumn{7}{>{\centering\arraybackslash}m{\cellWidtha + \cellWidthb + \cellWidthc + \cellWidthd + \cellWidthe + \cellWidthf + \cellWidthg+12\tabcolsep}}{\textbf{$\textbf{a}\ \textbf{=}\ \textbf{b}\ \textbf{=}\ \textbf{0.5}\textbf{,}\ \textbf{f}\ \textbf{=}\ \textbf{0.5}\textbf{,}\ \bm{\upmu}\ \textbf{=}\ \textbf{0.1}\textbf{,}\ \bm{\upomega}\ \textbf{=}\ \textbf{0.79} $}}\\
\cmidrule{1-7}
\textbf{T} & \textbf{$\bm{\upgamma}\ \textbf{=}\ \textbf{0.1} $} & \textbf{$\bm{\upgamma}\ \textbf{=}\ \textbf{0.2} $} & \textbf{$\bm{\upgamma}\ \textbf{=}\ \textbf{0.3} $} & \textbf{$\bm{\upgamma}\ \textbf{=}\ \textbf{0.5} $} & \textbf{$\bm{\upgamma}\ \textbf{=}\ \textbf{0.9} $} & \textbf{$\bm{\upgamma}\ \textbf{=}\ \textbf{1.0} $}\\
\cmidrule{1-7}
0.1 & $1.6 \times 10^{- 2} $ & $1.4 \times 10^{- 2} $ & $4.1 \times 10^{- 2} $ & $1.1 \times 10^{- 1} $ & $4.0 \times 10^{- 1} $ & $4.9 \times 10^{- 1} $\\
\cmidrule{1-7}
0.3 & $6.2 \times 10^{- 3} $ & $1.7 \times 10^{- 3} $ & $1.1 \times 10^{- 2} $ & $5.6 \times 10^{- 2} $ & $4.2 \times 10^{- 1} $ & $6.1 \times 10^{- 1} $\\
\cmidrule{1-7}
0.5 & $2.9 \times 10^{- 3} $ & $3.0 \times 10^{- 4} $ & $6.3 \times 10^{- 3} $ & $4.0 \times 10^{- 2} $ & $4.8 \times 10^{- 1} $ & $7.9 \times 10^{- 1} $\\
\cmidrule{1-7}
0.7 & $4.4 \times 10^{- 4} $ & $1.7 \times 10^{- 5} $ & $2.4 \times 10^{- 3} $ & $1.8 \times 10^{- 2} $ & $3.0 \times 10^{- 1} $ & $5.5 \times 10^{- 1} $\\
\cmidrule{1-7}
0.9 & $4.4 \times 10^{- 4} $ & $6.3 \times 10^{- 5} $ & $1.2 \times 10^{- 3} $ & $1.0 \times 10^{- 2} $ & $2.3 \times 10^{- 1} $ & $4.4 \times 10^{- 1} $\\
\bottomrule
\end{tabularx}}
\end{table}
\vspace{-12pt}
    
\begin{table}[H]
\tablesize{\small}
\caption{Estimated residual errors of solutions for different $\text{$\upalpha$}(\text{t}) $, k = 1 and $\text{$\upgamma$} = 0.2 $ in Example 1.}
\label{tabref:mathematics-2384216-t005}
\begin{adjustwidth}{-\extralength}{0cm}
\setlength{\cellWidtha}{\fulllength/9-2\tabcolsep-0in}
\setlength{\cellWidthb}{\fulllength/9-2\tabcolsep-0in}
\setlength{\cellWidthc}{\fulllength/9-2\tabcolsep-0in}
\setlength{\cellWidthd}{\fulllength/9-2\tabcolsep-0in}
\setlength{\cellWidthe}{\fulllength/9-2\tabcolsep-0in}
\setlength{\cellWidthf}{\fulllength/9-2\tabcolsep-0in}
\setlength{\cellWidthg}{\fulllength/9-2\tabcolsep-0in}
\setlength{\cellWidthh}{\fulllength/9-2\tabcolsep-0in}
\setlength{\cellWidthi}{\fulllength/9-2\tabcolsep-0in}
\scalebox{1}[1]{\begin{tabularx}{\fulllength}{>{\centering\arraybackslash}m{\cellWidtha}>{\centering\arraybackslash}m{\cellWidthb}>{\centering\arraybackslash}m{\cellWidthc}>{\centering\arraybackslash}m{\cellWidthd}>{\centering\arraybackslash}m{\cellWidthe}>{\centering\arraybackslash}m{\cellWidthf}>{\centering\arraybackslash}m{\cellWidthg}>{\centering\arraybackslash}m{\cellWidthh}>{\centering\arraybackslash}m{\cellWidthi}}
\toprule
\multicolumn{9}{>{\centering\arraybackslash}m{\cellWidtha + \cellWidthb + \cellWidthc + \cellWidthd + \cellWidthe + \cellWidthf + \cellWidthg + \cellWidthh + \cellWidthi+16\tabcolsep}}{\textbf{
Residual Errors, $\bm{\upgamma}\ \textbf{=}\ \textbf{0.2} $}}\\
\cmidrule{1-9}
\multirow{2}{*}{\parbox{\cellWidtha}{\centering \textbf{t}}\vspace{-4pt}} & \multicolumn{2}{>{\centering\arraybackslash}m{\cellWidthb + \cellWidthc+2\tabcolsep}}{\textbf{$\bm{\upalpha}(\textbf{t})\ \textbf{=}\ \textbf{1.2} $}} & \multicolumn{2}{>{\centering\arraybackslash}m{\cellWidthd + \cellWidthe+2\tabcolsep}}{\textbf{$\bm{\upalpha}(\textbf{t})\ \textbf{=}\ \textbf{1.4} $}} & \multicolumn{2}{>{\centering\arraybackslash}m{\cellWidthf + \cellWidthg+2\tabcolsep}}{\textbf{$\bm{\upalpha}(\textbf{t})\ \textbf{=}\ \textbf{1.6} $}} & \multicolumn{2}{>{\centering\arraybackslash}m{\cellWidthh + \cellWidthi+2\tabcolsep}}{\textbf{$\bm{\upalpha}(\textbf{t})\ \textbf{=}\ \textbf{1.8} $}}\\
\cmidrule{2-9}
 & \textbf{M = 3} & \textbf{M = 5} & \textbf{M = 3} & \textbf{M = 5} & \textbf{M = 3} & \textbf{M = 5} & \textbf{M = 3} & \textbf{M = 5}\\
\cmidrule{1-9}
0.1 & $2.7 \times 10^{- 1} $ & $8.3 \times 10^{- 3} $ & $1.8 \times 10^{- 2} $ & $1.5 \times 10^{- 3} $ & $2.7 \times 10^{- 1} $ & $2.5 \times 10^{- 4} $ & $6.4 \times 10^{- 1} $ & $1.4 \times 10^{- 3} $\\
\cmidrule{1-9}
0.3 & $6.0 \times 10^{- 3} $ & $1.5 \times 10^{- 3} $ & $5.2 \times 10^{- 4} $ & $1.3 \times 10^{- 4} $ & $5.5 \times 10^{- 3} $ & $1.6 \times 10^{- 4} $ & $9.7 \times 10^{- 3} $ & $1.1 \times 10^{- 4} $\\
\cmidrule{1-9}
0.5 & $2.4 \times 10^{- 2} $ & $7.3 \times 10^{- 4} $ & $8.8 \times 10^{- 3} $ & $4.7 \times 10^{- 5} $ & $2.9 \times 10^{- 2} $ & $1.0 \times 10^{- 4} $ & $4.4 \times 10^{- 2} $ & $1.8 \times 10^{- 5} $\\
\cmidrule{1-9}
0.7 & $2.6 \times 10^{- 4} $ & $2.8 \times 10^{- 4} $ & $6.2 \times 10^{- 4} $ & $1.6 \times 10^{- 5} $ & $1.1 \times 10^{- 3} $ & $4.4 \times 10^{- 5} $ & $1.4 \times 10^{- 3} $ & $3.2 \times 10^{- 6} $\\
\cmidrule{1-9}
0.9 & $4.6 \times 10^{- 3} $ & $1.5 \times 10^{- 4} $ & $6.2 \times 10^{- 3} $ & $9.0 \times 10^{- 6} $ & $7.0 \times 10^{- 3} $ & $2.6 \times 10^{- 5} $ & $7.4 \times 10^{- 3} $ & $6.2 \times 10^{- 6} $\\
\bottomrule
\end{tabularx}}
\end{adjustwidth}
\end{table}
\vspace{-12pt}
    
\begin{enumerate}
\item[(B)] \textbf{\boldmath{For Double-well}} $(\text{a} < 0,\ \text{b} > 0)$.
\end{enumerate}

\begin{figure}[H]
\includegraphics[scale=1.03]{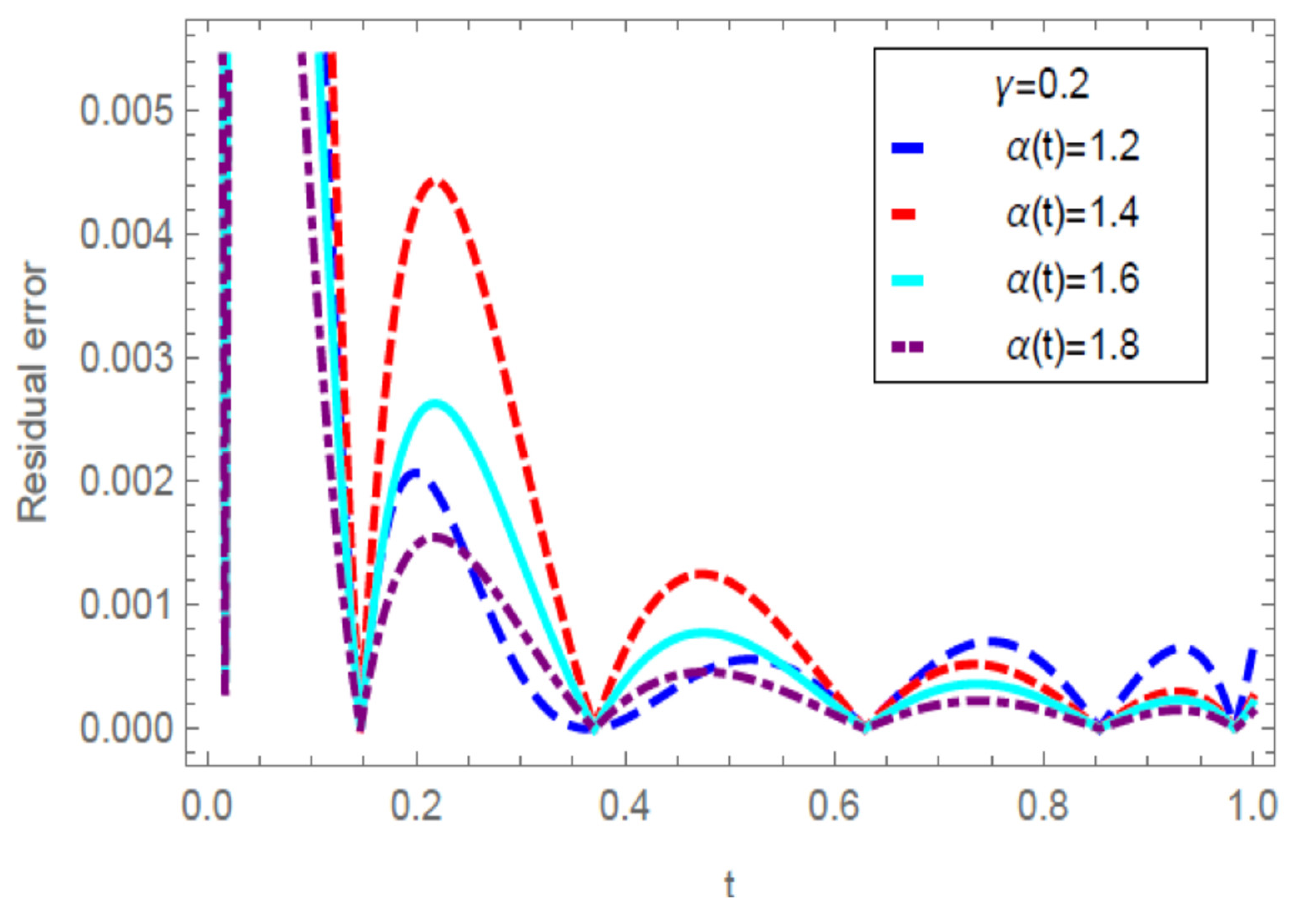}
\caption{Graph of Residual errors for double well case with k = 1, M = 5.}
\label{fig:mathematics-2384216-f004}
\end{figure}

The estimated AEs in the solutions of $\Im(\text{t}) $ with the comparison of the Legendre wavelet-Picard scheme (LWPS) and the ultraspherical wavelets scheme (UWS) for $\text{$\upalpha$}(\text{t}) = 2 $ and different FOBWs are listed in Table~\ref{tabref:mathematics-2384216-t006}. It can be easily analyzed from Table~\ref{tabref:mathematics-2384216-t006} that the suitable value of $\text{$\upgamma$} $ is 1 for achieving the best accuracy in the solution of the given model, with $\text{$\upalpha$}(\text{t}) = 2 $ and the proposed approach is superior to UWS~\cite{B28-mathematics-2384216} and LWPS~\cite{B29-mathematics-2384216} by considering the RK-4 solution~\cite{B28-mathematics-2384216} as an approximated analytical solution. The residual errors in the solutions of $\Im(\text{t}) $ for $\text{$\upalpha$}(\text{t}) = 1.5 $ and $\text{$\upalpha$}(\text{t}) = 1 + \sin \text{t} $ under different parameters mentioned are shown 
in Tables~\ref{tabref:mathematics-2384216-t007} and \ref{tabref:mathematics-2384216-t008}, respectively. In addition, the estimated residual errors in the solutions for $\text{$\upgamma$} = 0.2 $ and selected $\text{$\upalpha$}(\text{t}) $ are listed in Table~\ref{tabref:mathematics-2384216-t009}. The graphical interpretation of residual errors of solutions for the double-well case with selected values of $\text{$\upalpha$}(\text{t}) $, and $\text{$\upgamma$} = 0.2 $ is shown in Figure~\ref{fig:mathematics-2384216-f004}. The computed solutions are obtained for the first time with the variable order of the introduced model in terms of residual errors.    
\begin{table}[H]
\tablesize{\small}
\caption{Estimated AEs of solutions for $\text{$\upalpha$}(\text{t}) = 2 $ and selected $\text{$\upgamma$} $ in Example 1.}
\label{tabref:mathematics-2384216-t006}
    
\setlength{\cellWidtha}{\textwidth/6-2\tabcolsep-0in}
\setlength{\cellWidthb}{\textwidth/6-2\tabcolsep-0in}
\setlength{\cellWidthc}{\textwidth/6-2\tabcolsep-0in}
\setlength{\cellWidthd}{\textwidth/6-2\tabcolsep-0in}
\setlength{\cellWidthe}{\textwidth/6-2\tabcolsep-0in}
\setlength{\cellWidthf}{\textwidth/6-2\tabcolsep-0in}
\scalebox{1}[1]{\begin{tabularx}{\textwidth}{>{\centering\arraybackslash}m{\cellWidtha}>{\centering\arraybackslash}m{\cellWidthb}>{\centering\arraybackslash}m{\cellWidthc}>{\centering\arraybackslash}m{\cellWidthd}>{\centering\arraybackslash}m{\cellWidthe}>{\centering\arraybackslash}m{\cellWidthf}}
\toprule

\multicolumn{6}{>{\centering\arraybackslash}m{\cellWidtha + \cellWidthb + \cellWidthc + \cellWidthd + \cellWidthe + \cellWidthf+10\tabcolsep}}{\textbf{$\textbf{a}\ \textbf{=}\ \textbf{$-$}\textbf{0.5}\textbf{,}\ \textbf{b}\ \textbf{=}\ \textbf{0.5}\textbf{,}\ \textbf{f}\ \textbf{=}\ \textbf{0.5}\textbf{,}\ \bm{\upmu}\ \textbf{=}\ \textbf{0.1}\textbf{,}\ \bm{\upomega}\ \textbf{=}\ \textbf{0.79} $}}\\
\cmidrule{1-6}
\multirow{2}{*}{\parbox{\cellWidtha}{\centering \textbf{T}}\vspace{-4pt}} & \multicolumn{3}{>{\centering\arraybackslash}m{\cellWidthb + \cellWidthc + \cellWidthd+4\tabcolsep}}{\textbf{Proposed Approach, $\textbf{k}\ \textbf{=}\ \textbf{1}\textbf{,}\ \textbf{M}\ \textbf{=}\ \textbf{5} $}} & \multicolumn{2}{>{\centering\arraybackslash}m{\cellWidthe + \cellWidthf+2\tabcolsep}}{\textbf{Reference Approach, $\textbf{M}\ \textbf{=}\ \textbf{6} $}}\\
\cmidrule{2-6}
& \textbf{$\bm{\upgamma}\ \textbf{=}\ \textbf{0.5} $} & \textbf{$\bm{\upgamma}\ \textbf{=}\ \textbf{0.9} $} & \textbf{$\bm{\upgamma}\ 
\textbf{=}\ \textbf{1} $} & \textbf{UWS {\cite{B28-mathematics-2384216}}} & \textbf{LWPS {\cite{B29-mathematics-2384216}}}\\
\cmidrule{1-6}
0.1 & $1.5 \times 10^{- 6} $ & $2.0 \times 10^{- 8} $ & $7.9 \times 10^{- 8} $ & $1.1 \times 10^{- 8} $ & $1.0 \times 10^{- 8} $\\
\cmidrule{1-6}
0.3 & $7.4 \times 10^{- 6} $ & $8.2 \times 10^{- 7} $ & $1.9 \times 10^{- 8} $ & $9.0 \times 10^{- 8} $ & $5.9 \times 10^{- 8} $\\
\cmidrule{1-6}
0.5 & $9.0 \times 10^{- 6} $ & $6.6 \times 10^{- 7} $ & $7.9 \times 10^{- 8} $ & $3.6 \times 10^{- 7} $ & $1.7 \times 10^{- 7} $\\
\cmidrule{1-6}
0.7 & $1.4 \times 10^{- 5} $ & $1.3 \times 10^{- 6} $ & $7.0 \times 10^{- 8} $ & $5.6 \times 10^{- 7} $ & $3.6 \times 10^{- 7} $\\
\cmidrule{1-6}
0.9 & $1.7 \times 10^{- 5} $ & $1.5 \times 10^{- 6} $ & $7.2 \times 10^{- 8} $ & $1.1 \times 10^{- 6} $ & $9.4 \times 10^{- 7} $\\
\bottomrule
\end{tabularx}}
\end{table}
\vspace{-12pt}
    
\begin{table}[H]
\tablesize{\small}
\caption{\textls[-15]{Estimated residual errors of solutions for $\text{$\upalpha$}(\text{t}) = 1.5, $ $\text{k} = 1,\text{M} = 5, $ and selected $\text{$\upgamma$} $ in Example 1.}}
\label{tabref:mathematics-2384216-t007}
\setlength{\cellWidtha}{\textwidth/7-2\tabcolsep-0in}
\setlength{\cellWidthb}{\textwidth/7-2\tabcolsep-0in}
\setlength{\cellWidthc}{\textwidth/7-2\tabcolsep-0in}
\setlength{\cellWidthd}{\textwidth/7-2\tabcolsep-0in}
\setlength{\cellWidthe}{\textwidth/7-2\tabcolsep-0in}
\setlength{\cellWidthf}{\textwidth/7-2\tabcolsep-0in}
\setlength{\cellWidthg}{\textwidth/7-2\tabcolsep-0in}
\scalebox{1}[1]{
\begin{tabularx}{\textwidth}{>{\centering\arraybackslash}m{\cellWidtha}>{\centering\arraybackslash}m{\cellWidthb}>{\centering\arraybackslash}m{\cellWidthc}>{\centering\arraybackslash}m{\cellWidthd}>{\centering\arraybackslash}m{\cellWidthe}>{\centering\arraybackslash}m{\cellWidthf}>{\centering\arraybackslash}m{\cellWidthg}}
\toprule
\multicolumn{7}{>{\centering\arraybackslash}m{\cellWidtha + \cellWidthb + \cellWidthc + \cellWidthd + \cellWidthe + \cellWidthf + \cellWidthg+12\tabcolsep}}{\textbf{$\textbf{a}\ \textbf{=}\ \textbf{$-$}\textbf{0.5}\textbf{,}\ \textbf{b}\ \textbf{=}\ \textbf{0.5}\textbf{,}\ \textbf{f}\ \textbf{=}\ \textbf{0.5}\textbf{,}\ \bm{\upmu}\ \textbf{=}\ \textbf{0.1}\textbf{,}\ \bm{\upomega}\ \textbf{=}\ \textbf{0.79} $}}\\
\cmidrule{1-7}
\textbf{t} & \textbf{$\bm{\upgamma}\ \textbf{=}\ \textbf{0.1} $} & \textbf{$\bm{\upgamma}\ \textbf{=}\ \textbf{0.2} $} & \textbf{$\bm{\upgamma}\ \textbf{=}\ \textbf{0.3} $} & \textbf{$\bm{\upgamma}\ \textbf{=}\ \textbf{0.5} $} & \textbf{$\bm{\upgamma}\ \textbf{=}\ \textbf{0.9} $} & \textbf{$\bm{\upgamma}\ \textbf{=}\ \textbf{1.0} $}\\
\cmidrule{1-7}
0.1 & $1.2 \times 10^{- 2} $ & $8.7 \times 10^{- 3} $ & $7.7 \times 10^{- 3} $ & $1.6 \times 10^{- 2} $ & $6.9 \times 10^{- 2} $ & $8.8 \times 10^{- 2} $\\
\cmidrule{1-7}
0.3 & $2.2 \times 10^{- 3} $ & $1.8 \times 10^{- 3} $ & $1.6 \times 10^{- 3} $ & $4.0 \times 10^{- 3} $ & $3.6 \times 10^{- 2} $ & $5.6 \times 10^{- 2} $\\
\cmidrule{1-7}
0.5 & $1.1 \times 10^{- 3} $ & $9.6 \times 10^{- 4} $ & $8.8 \times 10^{- 4} $ & $2.3 \times 10^{- 3} $ & $3.3 \times 10^{- 2} $ & $5.6 \times 10^{- 2} $\\
\cmidrule{1-7}
0.7 & $4.6 \times 10^{- 4} $ & $4.0 \times 10^{- 4} $ & $3.7 \times 10^{- 4} $ & $1.0 \times 10^{- 3} $ & $2.0 \times 10^{- 2} $ & $3.7 \times 10^{- 2} $\\
\cmidrule{1-7}
0.9 & $2.6 \times 10^{- 4} $ & $2.3 \times 10^{- 4} $ & $2.1 \times 10^{- 4} $ & $6.2 \times 10^{- 4} $ & $1.6 \times 10^{- 2} $ & $3.2 \times 10^{- 2} $\\
\bottomrule
\end{tabularx}}
\end{table}
\vspace{-12pt}
    
\begin{table}[H]
\tablesize{\small}
\caption{Estimated residual errors of solutions for $\text{$\upalpha$}(\text{t}) 
= 1 + \sin \text{t}, $ $\text{k} = 1,\text{M} = 5, $ and selected $\text{$\upgamma$} $ in Example 1.}
\label{tabref:mathematics-2384216-t008}
\setlength{\cellWidtha}{\textwidth/7-2\tabcolsep-0in}
\setlength{\cellWidthb}{\textwidth/7-2\tabcolsep-0in}
\setlength{\cellWidthc}{\textwidth/7-2\tabcolsep-0in}
\setlength{\cellWidthd}{\textwidth/7-2\tabcolsep-0in}
\setlength{\cellWidthe}{\textwidth/7-2\tabcolsep-0in}
\setlength{\cellWidthf}{\textwidth/7-2\tabcolsep-0in}
\setlength{\cellWidthg}{\textwidth/7-2\tabcolsep-0in}
\scalebox{1}[1]{\begin{tabularx}{\textwidth}{>{\centering\arraybackslash}m{\cellWidtha}>{\centering\arraybackslash}m{\cellWidthb}>{\centering\arraybackslash}m{\cellWidthc}>{\centering\arraybackslash}m{\cellWidthd}>{\centering\arraybackslash}m{\cellWidthe}>{\centering\arraybackslash}m{\cellWidthf}>{\centering\arraybackslash}m{\cellWidthg}}
\toprule
\multicolumn{7}{>{\centering\arraybackslash}m{\cellWidtha + \cellWidthb + \cellWidthc + \cellWidthd + \cellWidthe + \cellWidthf + \cellWidthg+12\tabcolsep}}{\textbf{$\textbf{a}\ \textbf{=}\ \textbf{$-$}\textbf{0.5}\textbf{,}\ \textbf{b}\ \textbf{=}\ \textbf{0.5}\textbf{,}\ \textbf{f}\ \textbf{=}\ \textbf{0.5}\textbf{,}\ \bm{\upmu}\ \textbf{=}\ \textbf{0.1}\textbf{,}\ \bm{\upomega}\ \textbf{=}\ \textbf{0.79} $}}\\
\cmidrule{1-7}
\textbf{t} & \textbf{$\bm{\upgamma}\ \textbf{=}\ \textbf{0.1} $} & \textbf{$\bm{\upgamma}\ \textbf{=}\ \textbf{0.2} $} & \textbf{$\bm{\upgamma}\ \textbf{=}\ \textbf{0.3} $} & \textbf{$\bm{\upgamma}\ \textbf{=}\ \textbf{0.5} $} & \textbf{$\bm{\upgamma}\ \textbf{=}\ \textbf{0.9} $} & \textbf{$\bm{\upgamma}\ \textbf{=}\ \textbf{1.0} $}\\
\cmidrule{1-7}
0.1 & $5.3 \times 10^{- 2} $ & $4.2 \times 10^{- 2} $ & $5.2 \times 10^{- 2} $ & $1.1 \times 10^{- 1} $ & $3.9 \times 10^{- 1} $ & $4.7 \times 10^{- 1} $\\
\cmidrule{1-7}
0.3 & $1.7 \times 10^{- 2} $ & $1.3 \times 10^{- 2} $ & $1.7 \times 10^{- 2} $ & $5.5 \times 10^{- 2} $ & $4.0 \times 10^{- 1} $ & $5.9 \times 10^{- 1} $\\
\cmidrule{1-7}
0.5 & $1.0 \times 10^{- 2} $ & $7.8 \times 10^{- 3} $ & $1.0 \times 10^{- 2} $ & $3.9 \times 10^{- 2} $ & $4.5 \times 10^{- 1} $ & $7.6 \times 10^{- 1} $\\
\cmidrule{1-7}
0.7 & $4.0 \times 10^{- 3} $ & $2.8 \times 10^{- 3} $ & $3.9 \times 10^{- 3} $ & $1.7 \times 10^{- 2} $ & $2.9 \times 10^{- 1} $ & $5.4 \times 10^{- 1} $\\
\cmidrule{1-7}
0.9 & $1.8 \times 10^{- 3} $ & $1.0 \times 10^{- 3} $ & $1.7 \times 10^{- 3} $ & $1.0 \times 10^{- 2} $ & $2.2 \times 10^{- 1} $ & $4.4 \times 10^{- 1} $\\
\bottomrule
\end{tabularx}}
\end{table}
\vspace{-6pt}
    
\begin{table}[H]
\tablesize{\small}
\caption{Estimated residual errors of solutions for different $\text{$\upalpha$}(\text{t})$, 
$k = 1$ and $\text{$\upgamma$} = 0.2 $ in Example 1.}
\label{tabref:mathematics-2384216-t009}
    
\begin{adjustwidth}{-\extralength}{0cm}
\setlength{\cellWidtha}{\fulllength/9-2\tabcolsep-0in}
\setlength{\cellWidthb}{\fulllength/9-2\tabcolsep-0in}
\setlength{\cellWidthc}{\fulllength/9-2\tabcolsep-0in}
\setlength{\cellWidthd}{\fulllength/9-2\tabcolsep-0in}
\setlength{\cellWidthe}{\fulllength/9-2\tabcolsep-0in}
\setlength{\cellWidthf}{\fulllength/9-2\tabcolsep-0in}
\setlength{\cellWidthg}{\fulllength/9-2\tabcolsep-0in}
\setlength{\cellWidthh}{\fulllength/9-2\tabcolsep-0in}
\setlength{\cellWidthi}{\fulllength/9-2\tabcolsep-0in}
\scalebox{1}[1]{\begin{tabularx}{\fulllength}{>{\centering\arraybackslash}m{\cellWidtha}>{\centering\arraybackslash}m{\cellWidthb}>{\centering\arraybackslash}m{\cellWidthc}>{\centering\arraybackslash}m{\cellWidthd}>{\centering\arraybackslash}m{\cellWidthe}>{\centering\arraybackslash}m{\cellWidthf}>{\centering\arraybackslash}m{\cellWidthg}>{\centering\arraybackslash}m{\cellWidthh}>{\centering\arraybackslash}m{\cellWidthi}}
\toprule
\multicolumn{9}{>{\centering\arraybackslash}m{\cellWidtha + \cellWidthb + \cellWidthc + \cellWidthd + \cellWidthe + \cellWidthf + \cellWidthg + \cellWidthh + \cellWidthi+16\tabcolsep}}{\textbf{
Residual Errors, $\bm{\upgamma}\ \textbf{=}\ \textbf{0.2} $}}\\
\cmidrule{1-9}
\multirow{2}{*}{\parbox{\cellWidtha}{\centering \textbf{t}}\vspace{-4pt}} & \multicolumn{2}{>{\centering\arraybackslash}m{\cellWidthb + \cellWidthc+2\tabcolsep}}{\textbf{$\bm{\upalpha}(\textbf{t})\ \textbf{=}\ \textbf{1.2} $}} & \multicolumn{2}{>{\centering\arraybackslash}m{\cellWidthd + \cellWidthe+2\tabcolsep}}{\textbf{$\bm{\upalpha}(\textbf{t})\ \textbf{=}\ \textbf{1.4} $}} & \multicolumn{2}{>{\centering\arraybackslash}m{\cellWidthf + \cellWidthg+2\tabcolsep}}{\textbf{$\bm{\upalpha}(\textbf{t})\ \textbf{=}\ \textbf{1.6} $}} & \multicolumn{2}{>{\centering\arraybackslash}m{\cellWidthh + \cellWidthi+2\tabcolsep}}{\textbf{$\bm{\upalpha}(\textbf{t})\ \textbf{=}\ \textbf{1.8} $}}\\
\cmidrule{2-9}
& \textbf{M = 3} & \textbf{M = 5} & \textbf{M = 3} & \textbf{M = 5} & \textbf{M = 3} & \textbf{M = 5} & \textbf{M = 3} & \textbf{M = 5}\\
\cmidrule{1-9}
0.1 & $3.2 \times 10^{- 1} $ & $8.6 \times 10^{- 3} $ & $1.2 \times 10^{- 1} $ & $1.1 \times 10^{- 2} $ & $1.0 \times 10^{- 1} $ & $6.7 \times 10^{- 3} $ & $4.0 \times 10^{- 1} $ & $4.1 \times 10^{- 3} $\\
\cmidrule{1-9}
0.3 & $3.5 \times 10^{- 3} $ & $4.3 \times 10^{- 4} $ & $2.4 \times 10^{- 3} $ & $2.3 \times 10^{- 3} $ & $7.5 \times 10^{- 3} $ & $1.3 \times 10^{- 3} $ & $1.2 \times 10^{- 2} $ & $8.0 \times 10^{- 4} $\\
\cmidrule{1-9}
0.5 & $1.7 \times 10^{- 2} $ & $5.4 \times 10^{- 4} $ & $4.6 \times 10^{- 2} $ & $1.1 \times 10^{- 3} $ & $6.9 \times 10^{- 2} $ & $7.4 \times 10^{- 4} $ & $8.9 \times 10^{- 2} $ & $4.3 \times 10^{- 4} $\\
\cmidrule{1-9}
0.7 & $2.7 \times 10^{- 3} $ & $5.5 \times 10^{- 4} $ & $3.3 \times 10^{- 3} $ & $4.6 \times 10^{- 4} $ & $3.8 \times 10^{- 3} $ & $3.1 \times 10^{- 4} $ & $4.2 \times 10^{- 3} $ & $1.9 \times 10^{- 4} $\\
\cmidrule{1-9}
0.9 & $3.5 \times 10^{- 2} $ & $4.9 \times 10^{- 4} $ & $3.3 \times 10^{- 2} $ & $2.4 \times 10^{- 4} $ & $3.2 \times 10^{- 2} $ & $1.8 \times 10^{- 4} $ & $3.2 \times 10^{- 2} $ & $1.2 \times 10^{- 4} $\\
\bottomrule
\end{tabularx}}
\end{adjustwidth}
\end{table}
\vspace{-12pt}
    
\begin{enumerate}
\item[(C)] \textbf{\boldmath{For Double-hump well}} $(\text{a} > 0,\ \text{b} < 0)$.
\end{enumerate}

\begin{figure}[H]
\includegraphics[scale=1]{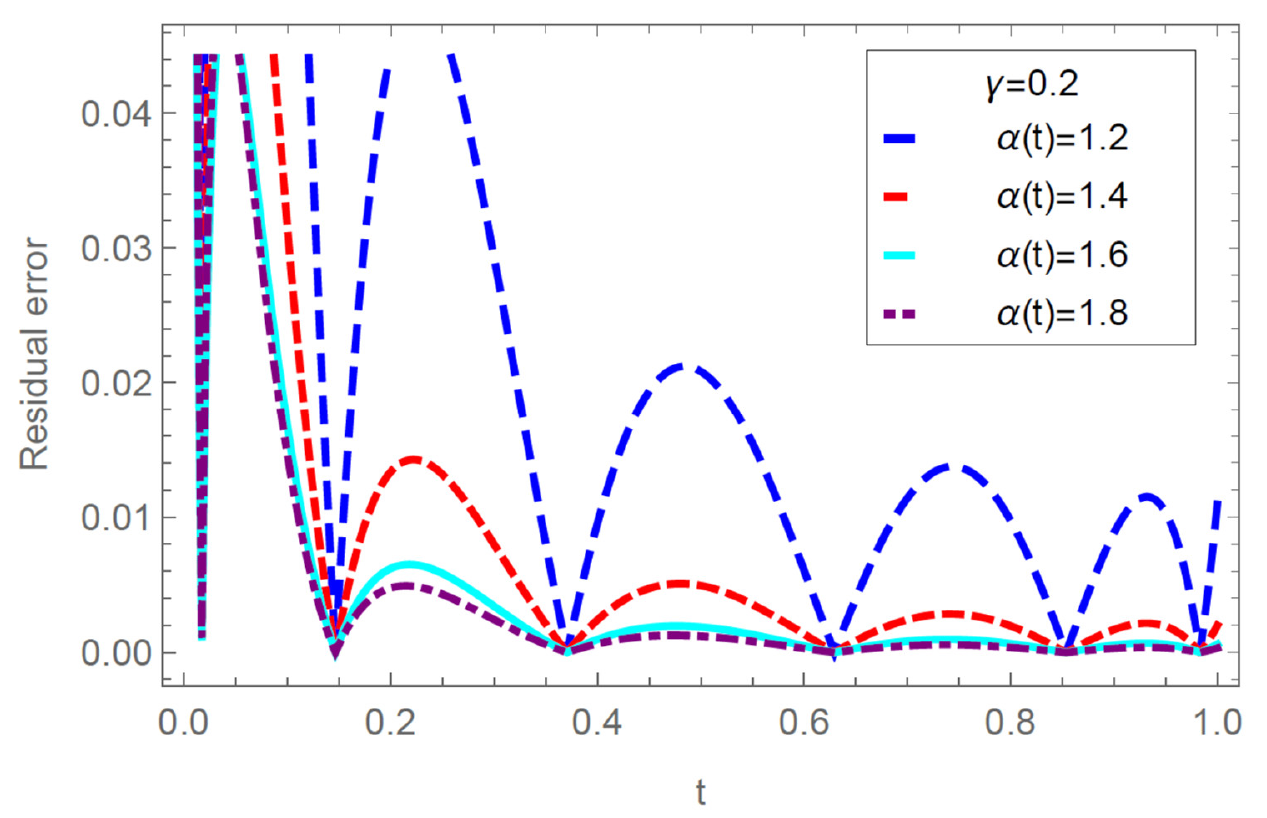}
\caption{Graph of Residual errors for double-hump case with k = 1, M = 5.}
\label{fig:mathematics-2384216-f005}
\end{figure}

The estimated AEs in the solutions of $\Im(\text{t}) $ with the comparison of the Legendre wavelet-Picard scheme (LWPS) and the ultraspherical wavelets scheme (UWS) for $\text{$\upalpha$}(\text{t}) = 2 $ and different FOBW bases are listed in Table~\ref{tabref:mathematics-2384216-t010}. It can be easily analyzed from Table~\ref{tabref:mathematics-2384216-t010} that the suitable value of $\text{$\upgamma$} $ is 1 for achieving the best accuracy in the solution of the given model, with $\text{$\upalpha$}(\text{t}) = 2 $ and the proposed approach is superior to UWS~\cite{B28-mathematics-2384216} and LWPS~\cite{B29-mathematics-2384216} by considering the RK-4 solution~\cite{B28-mathematics-2384216} as an approximated analytical solution. The residual error in the solutions of $\Im(\text{t}) $ for $\text{$\upalpha$}(\text{t}) = 1.5 $ and $\text{$\upalpha$}(\text{t}) = 1 + \sin \text{t} $ with selected parameters is presented 
in Tables~\ref{tabref:mathematics-2384216-t011} and \ref{tabref:mathematics-2384216-t012}, respectively. Furthermore, the estimated residual errors in the solutions for $\text{$\upgamma$} = 0.2 $ and selected $\text{$\upalpha$}(\text{t}) $ are given in Table~\ref{tabref:mathematics-2384216-t013}. The graphical interpretation of residual errors of solutions for the double hump case with selected values of $\text{$\upalpha$}(\text{t}) $, and $\text{$\upgamma$} = 0.2 $ is shown in Figure~\ref{fig:mathematics-2384216-f005}. The computed solutions are obtained for the first time with the variable order of the introduced model in terms of residual errors.
\begin{table}[H]
\tablesize{\small}
\caption{Estimated AEs of solutions for $\text{$\upalpha$}(\text{t}) = 2 $ and selected $\text{$\upgamma$} $ in Example 1.}
\label{tabref:mathematics-2384216-t010}
\setlength{\cellWidtha}{\textwidth/6-2\tabcolsep-0in}
\setlength{\cellWidthb}{\textwidth/6-2\tabcolsep-0in}
\setlength{\cellWidthc}{\textwidth/6-2\tabcolsep-0in}
\setlength{\cellWidthd}{\textwidth/6-2\tabcolsep-0in}
\setlength{\cellWidthe}{\textwidth/6-2\tabcolsep-0in}
\setlength{\cellWidthf}{\textwidth/6-2\tabcolsep-0in}
    
\scalebox{1}[1]{
\begin{tabularx}{\textwidth}{>{\centering\arraybackslash}m{\cellWidtha}>{\centering\arraybackslash}m{\cellWidthb}>{\centering\arraybackslash}m{\cellWidthc}>{\centering\arraybackslash}m{\cellWidthd}>{\centering\arraybackslash}m{\cellWidthe}>{\centering\arraybackslash}m{\cellWidthf}}
\toprule
\multicolumn{6}{>{\centering\arraybackslash}m{\cellWidtha + \cellWidthb + \cellWidthc + \cellWidthd + \cellWidthe + \cellWidthf+10\tabcolsep}}{\textbf{$\textbf{a}\ \textbf{=}\ \textbf{0.5}\textbf{,}\ \textbf{b}\ \textbf{=}\ \textbf{$-$}\textbf{0.5}\textbf{,}\ \textbf{f}\ \textbf{=}\ \textbf{0.5}\textbf{,}\ \bm{\upmu}\ \textbf{=}\ \textbf{0.1}\textbf{,}\ \bm{\upomega}\ \textbf{=}\ \textbf{0.79} $}}\\
\cmidrule{1-6}
\multirow{2}{*}{\parbox{\cellWidtha}{\centering \textbf{t}}\vspace{-4pt}} & \multicolumn{3}{>{\centering\arraybackslash}m{\cellWidthb + \cellWidthc + \cellWidthd+4\tabcolsep}}{\textbf{Proposed Approach, $\textbf{k}\ \textbf{=}\ \textbf{1}\textbf{,}\ \textbf{M}\ \textbf{=}\ \textbf{5} $}} & \multicolumn{2}{>{\centering\arraybackslash}m{\cellWidthe + \cellWidthf+2\tabcolsep}}{\textbf{Reference Approach, $\textbf{M}\ \textbf{=}\ \textbf{6} $}}\\
\cmidrule{2-6}
& \textbf{$\bm{\upgamma}\ \textbf{=}\ \textbf{0.5} $} & \textbf{$\bm{\upgamma}\ \textbf{=}\ \textbf{0.9} $} & \textbf{$\bm{\upgamma}\ \textbf{=}\ \textbf{1} $} & \textbf{UWS {\cite{B28-mathematics-2384216}}
} & \textbf{LWPS {\cite{B29-mathematics-2384216}}
}\\
\cmidrule{1-6}
0.1 & $3.7 \times 10^{- 6} $ & $3.7 \times 10^{- 8} $ & $2.0 \times 10^{- 8} $ & $1.0 \times 10^{- 7} $ & $8.0 \times 10^{- 8} $\\
\cmidrule{1-6}
0.3 & $1.3 \times 10^{- 5} $ & $1.8 \times 10^{- 7} $ & $8.4 \times 10^{- 8} $ & $4.8 \times 10^{- 7} $ & $4.1 \times 10^{- 7} $\\
\cmidrule{1-6}
0.5 & $1.7 \times 10^{- 5} $ & $6.4 \times 10^{- 8} $ & $9.8 \times 10^{- 8} $ & $1.4 \times 10^{- 6} $ & $1.2 \times 10^{- 6} $\\
\cmidrule{1-6}
0.7 & $2.9 \times 10^{- 5} $ & $3.4 \times 10^{- 7} $ & $1.8 \times 10^{- 7} $ & $3.8 \times 10^{- 6} $ & $3.6 \times 10^{- 6} $\\
\cmidrule{1-6}
0.9 & $4.1 \times 10^{- 5} $ & $2.7 \times 10^{- 7} $ & $9.1 \times 10^{- 8} $ & $5.8 \times 10^{- 6} $ & $5.1 \times 10^{- 6} $\\
\bottomrule
\end{tabularx}}
\end{table}
\vspace{-12pt}

\begin{table}[H]
\tablesize{\small}
\caption{Estimated residual errors of solutions for $\text{$\upalpha$}(\text{t}) = 1.5, $ $\text{k} = 1,\text{M} = 5, $ and selected $\text{$\upgamma$} $ in Example 1.}
\label{tabref:mathematics-2384216-t011}
\setlength{\cellWidtha}{\textwidth/7-2\tabcolsep-0in}
\setlength{\cellWidthb}{\textwidth/7-2\tabcolsep-0in}
\setlength{\cellWidthc}{\textwidth/7-2\tabcolsep-0in}
\setlength{\cellWidthd}{\textwidth/7-2\tabcolsep-0in}
\setlength{\cellWidthe}{\textwidth/7-2\tabcolsep-0in}
\setlength{\cellWidthf}{\textwidth/7-2\tabcolsep-0in}
\setlength{\cellWidthg}{\textwidth/7-2\tabcolsep-0in}
    
\scalebox{1}[1]{
\begin{tabularx}{\textwidth}{>{\centering\arraybackslash}m{\cellWidtha}>{\centering\arraybackslash}m{\cellWidthb}>{\centering\arraybackslash}m{\cellWidthc}>{\centering\arraybackslash}m{\cellWidthd}>{\centering\arraybackslash}m{\cellWidthe}>{\centering\arraybackslash}m{\cellWidthf}>{\centering\arraybackslash}m{\cellWidthg}}
\toprule
\multicolumn{7}{>{\centering\arraybackslash}m{\cellWidtha + \cellWidthb + \cellWidthc + \cellWidthd + \cellWidthe + \cellWidthf + \cellWidthg+12\tabcolsep}}{\textbf{$\textbf{a}\ \textbf{=}\ \textbf{0.5}\textbf{,}\ \textbf{b}\ \textbf{=}\ \textbf{$-$}\textbf{0.5}\textbf{,}\ \textbf{f}\ \textbf{=}\ \textbf{0.5}\textbf{,}\ \bm{\upmu}\ \textbf{=}\ \textbf{0.1}\textbf{,}\ \bm{\upomega}\ \textbf{=}\ \textbf{0.79} $}}\\
\cmidrule{1-7}
\textbf{t} & \textbf{$\bm{\upgamma}\ \textbf{=}\ \textbf{0.1} $} & \textbf{$\bm{\upgamma}\ \textbf{=}\ \textbf{0.2} $} & \textbf{$\bm{\upgamma}\ \textbf{=}\ \textbf{0.3} $} & \textbf{$\bm{\upgamma}\ \textbf{=}\ \textbf{0.5} $} & \textbf{$\bm{\upgamma}\ \textbf{=}\ \textbf{0.9} $} & \textbf{$\bm{\upgamma}\ \textbf{=}\ \textbf{1.0} $}\\
\cmidrule{1-7}
0.1 & $2.8 \times 10^{- 2} $ & $2.1 \times 10^{- 2} $ & $1.3 \times 10^{- 2} $ & $1.7 \times 10^{- 2} $ & $7.0 \times 10^{- 2} $ & $8.9 \times 10^{- 2} $\\
\cmidrule{1-7}
0.3 & $4.9 \times 10^{- 3} $ & $4.8 \times 10^{- 3} $ & $3.6 \times 10^{- 3} $ & $4.8 \times 10^{- 3} $ & $3.7 \times 10^{- 2} $ & $5.7 \times 10^{- 2} $\\
\cmidrule{1-7}
0.5 & $2.5 \times 10^{- 3} $ & $2.8 \times 10^{- 3} $ & $2.3 \times 10^{- 3} $ & $3.0 \times 10^{- 3} $ & $3.4 \times 10^{- 2} $ & $5.8 \times 10^{- 2} $\\
\cmidrule{1-7}
0.7 & $1.0 \times 10^{- 3} $ & $1.3 \times 10^{- 3} $ & $1.1 \times 10^{- 3} $ & $1.5 \times 10^{- 3} $ & $2.1 \times 10^{- 2} $ & $3.8 \times 10^{- 2} $\\
\cmidrule{1-7}
0.9 & $6.7 \times 10^{- 4} $ & $8.8 \times 10^{- 4} $ & $8.1 \times 10^{- 4} $ & $1.0 \times 10^{- 3} $ & $1.6 \times 10^{- 2} $ & $3.3 \times 10^{- 2} $\\
\bottomrule
\end{tabularx}}
\end{table}
\vspace{-12pt}
    
\begin{table}[H]
\tablesize{\small}
\caption{Estimated residual errors for $\text{$\upalpha$}(\text{t}) = 1 + \sin \text{t}, $ $\text{k} = 1,\text{M} = 5, $ and selected $\text{$\upgamma$} $ in Example 1.}
\label{tabref:mathematics-2384216-t012}
\setlength{\cellWidtha}{\textwidth/7-2\tabcolsep-0in}
\setlength{\cellWidthb}{\textwidth/7-2\tabcolsep-0in}
\setlength{\cellWidthc}{\textwidth/7-2\tabcolsep-0in}
\setlength{\cellWidthd}{\textwidth/7-2\tabcolsep-0in}
\setlength{\cellWidthe}{\textwidth/7-2\tabcolsep-0in}
\setlength{\cellWidthf}{\textwidth/7-2\tabcolsep-0in}
\setlength{\cellWidthg}{\textwidth/7-2\tabcolsep-0in}
\scalebox{1}[1]{\begin{tabularx}{\textwidth}{>{\centering\arraybackslash}m{\cellWidtha}>{\centering\arraybackslash}m{\cellWidthb}>{\centering\arraybackslash}m{\cellWidthc}>{\centering\arraybackslash}m{\cellWidthd}>{\centering\arraybackslash}m{\cellWidthe}>{\centering\arraybackslash}m{\cellWidthf}>{\centering\arraybackslash}m{\cellWidthg}}
\toprule
\multicolumn{7}{>{\centering\arraybackslash}m{\cellWidtha + \cellWidthb + \cellWidthc + \cellWidthd + \cellWidthe + \cellWidthf + \cellWidthg+12\tabcolsep}}{\textbf{$\textbf{a}\ \textbf{=}\ \textbf{0.5}\textbf{,}\ \textbf{b}\ \textbf{=}\ \textbf{$-$}\textbf{0.5}\textbf{,}\ \textbf{f}\ \textbf{=}\ \textbf{0.5}\textbf{,}\ \bm{\upmu}\ \textbf{=}\ \textbf{0.1}\textbf{,}\ \bm{\upomega}\ \textbf{=}\ \textbf{0.79} $}}\\
\cmidrule{1-7}
\textbf{t} & \textbf{$\bm{\upgamma}\ \textbf{=}\ \textbf{0.1} $} & \textbf{$\bm{\upgamma}\ \textbf{=}\ \textbf{0.2} $} & \textbf{$\bm{\upgamma}\ \textbf{=}\ \textbf{0.3} $} & \textbf{$\bm{\upgamma}\ \textbf{=}\ \textbf{0.5} $} & \textbf{$\bm{\upgamma}\ \textbf{=}\ \textbf{0.9} $} & \textbf{$\bm{\upgamma}\ \textbf{=}\ \textbf{1.0} $}\\
\cmidrule{1-7}
0.1 & $1.3 \times 10^{- 1} $ & $2.5 \times 10^{- 1} $ & $3.2 \times 10^{- 1} $ & $2.2 \times 10^{- 1} $ & $3.9 \times 10^{- 1} $ & $4.8 \times 10^{- 1} $\\
\cmidrule{1-7}
0.3 & $5.2 \times 10^{- 2} $ & $1.2 \times 10^{- 1} $ & $2.0 \times 10^{- 1} $ & $1.8 \times 10^{- 1} $ & $4.2 \times 10^{- 1} $ & $6.1 \times 10^{- 1} $\\
\cmidrule{1-7}
0.5 & $3.6 \times 10^{- 2} $ & $1.0 \times 10^{- 1} $ & $1.9 \times 10^{- 1} $ & $2.0 \times 10^{- 1} $ & $5.1 \times 10^{- 1} $ & $8.1 \times 10^{- 1} $\\
\cmidrule{1-7}
0.7 & $1.7 \times 10^{- 2} $ & $5.5 \times 10^{- 2} $ & $1.1 \times 10^{- 1} $ & $1.4 \times 10^{- 1} $ & $3.6 \times 10^{- 1} $ & $6.0 \times 10^{- 1} $\\
\cmidrule{1-7}
0.9 & $1.1 \times 10^{- 2} $ & $4.0 \times 10^{- 2} $ & $9.6 \times 10^{- 2} $ & $1.3 \times 10^{- 1} $ & $3.1 \times 10^{- 1} $ & $5.2 \times 10^{- 1} $\\
\bottomrule
\end{tabularx}}
\end{table}
\vspace{-12pt}
   
\begin{table}[H]
\tablesize{\small}
\caption{Estimated residual errors of solutions for different $\text{$\upalpha$}(\text{t}) $, k = 1 and $\text{$\upgamma$} = 0.2 $ in Example 1.}
\label{tabref:mathematics-2384216-t013}
   
\begin{adjustwidth}{-\extralength}{0cm}
\setlength{\cellWidtha}{\fulllength/9-2\tabcolsep-0in}
\setlength{\cellWidthb}{\fulllength/9-2\tabcolsep-0in}
\setlength{\cellWidthc}{\fulllength/9-2\tabcolsep-0in}
\setlength{\cellWidthd}{\fulllength/9-2\tabcolsep-0in}
\setlength{\cellWidthe}{\fulllength/9-2\tabcolsep-0in}
\setlength{\cellWidthf}{\fulllength/9-2\tabcolsep-0in}
\setlength{\cellWidthg}{\fulllength/9-2\tabcolsep-0in}
\setlength{\cellWidthh}{\fulllength/9-2\tabcolsep-0in}
\setlength{\cellWidthi}{\fulllength/9-2\tabcolsep-0in}
       
\scalebox{1}[1]{\begin{tabularx}{\fulllength}{>{\centering\arraybackslash}m{\cellWidtha}>{\centering\arraybackslash}m{\cellWidthb}>{\centering\arraybackslash}m{\cellWidthc}>{\centering\arraybackslash}m{\cellWidthd}>{\centering\arraybackslash}m{\cellWidthe}>{\centering\arraybackslash}m{\cellWidthf}>{\centering\arraybackslash}m{\cellWidthg}>{\centering\arraybackslash}m{\cellWidthh}>{\centering\arraybackslash}m{\cellWidthi}}
\toprule
\multicolumn{9}{>{\centering\arraybackslash}m{\cellWidtha + \cellWidthb + \cellWidthc + \cellWidthd + \cellWidthe + \cellWidthf + \cellWidthg + \cellWidthh + \cellWidthi+16\tabcolsep}}{\textbf{
Residual Errors, $\bm{\upgamma}\ \textbf{=}\ \textbf{0.2} $}}\\
\cmidrule{1-9}
\multirow{2}{*}{\parbox{\cellWidtha}{\centering \textbf{t}}\vspace{-4pt}} & \multicolumn{2}{>{\centering\arraybackslash}m{\cellWidthb + \cellWidthc+2\tabcolsep}}{\textbf{$\bm{\upalpha}(\textbf{t})\ \textbf{=}\ \textbf{1.2} $}} & \multicolumn{2}{>{\centering\arraybackslash}m{\cellWidthd + \cellWidthe+2\tabcolsep}}{\textbf{$\bm{\upalpha}(\textbf{t})\ \textbf{=}\ \textbf{1.4} $}} & \multicolumn{2}{>{\centering\arraybackslash}m{\cellWidthf + \cellWidthg+2\tabcolsep}}{\textbf{$\bm{\upalpha}(\textbf{t})\ \textbf{=}\ \textbf{1.6} $}} & \multicolumn{2}{>{\centering\arraybackslash}m{\cellWidthh + \cellWidthi+2\tabcolsep}}{\textbf{$\bm{\upalpha}(\textbf{t})\ \textbf{=}\ \textbf{1.8} $}}\\
\cmidrule{2-9}
& \textbf{M = 3} & \textbf{M = 5} & \textbf{M = 3} & \textbf{M = 5} & \textbf{M = 3} & \textbf{M = 5} & \textbf{M = 3} & \textbf{M = 5}\\
\cmidrule{1-9}
0.1 & $1.0 \times 10^{+ 0} $ & $8.9 \times 10^{- 2} $ & $1.4 \times 10^{+ 0} $ & $3.1 \times 10^{- 2} $ & $2.1 \times 10^{+ 0} $ & $1.6 \times 10^{- 2} $ & $3.1 \times 10^{+ 0} $ & $1.4 \times 10^{- 2} $\\
\cmidrule{1-9}
0.3 & $1.2 \times 10^{- 2} $ & $2.9 \times 10^{- 2} $ & $1.2 \times 10^{- 2} $ & $7.9 \times 10^{- 3} $ & $1.3 \times 10^{- 2} $ & $3.4 \times 10^{- 3} $ & $1.3 \times 10^{- 2} $ & $2.4 \times 10^{- 3} $\\
\cmidrule{1-9}
0.5 & $4.5 \times 10^{- 2} $ & $2.0 \times 10^{- 2} $ & $7.1 \times 10^{- 2} $ & $4.9 \times 10^{- 3} $ & $9.3 \times 10^{- 2} $ & $1.8 \times 10^{- 3} $ & $1.1 \times 10^{- 1} $ & $1.1 \times 10^{- 3} $\\
\cmidrule{1-9}
0.7 & $8.0 \times 10^{- 2} $ & $1.1 \times 10^{- 2} $ & $9.4 \times 10^{- 3} $ & $2.4 \times 10^{- 3} $ & $1.0 \times 10^{- 2} $ & $8.3 \times 10^{- 4} $ & $1.1 \times 10^{- 2} $ & $4.9 \times 10^{- 4} $\\
\cmidrule{1-9}
0.9 & $1.0 \times 10^{- 1} $ & $8.8 \times 10^{- 3} $ & $1.1 \times 10^{- 1} $ & $1.6 \times 10^{- 3} $ & $1.2 \times 10^{- 1} $ & $5.4 \times 10^{- 4} $ & $1.2 \times 10^{- 1} $ & $3.0 \times 10^{- 4} $\\
\bottomrule
\end{tabularx}}
\end{adjustwidth}
\end{table}
\vspace{-12pt}
    
It can be concluded from Tables~\ref{tabref:mathematics-2384216-t003}, \ref{tabref:mathematics-2384216-t004},
\ref{tabref:mathematics-2384216-t005}, \ref{tabref:mathematics-2384216-t007},
\ref{tabref:mathematics-2384216-t008},\ref{tabref:mathematics-2384216-t009}, \ref{tabref:mathematics-2384216-t011},
\ref{tabref:mathematics-2384216-t012}, \ref{tabref:mathematics-2384216-t013} that the FOBW basis $(\text{$\upgamma$} \neq 1) $ 
is more efficient and reliable than the non-fractional wavelets $(\text{$\upgamma$} = 1) $ for solving this non-linear model.

\vspace{12pt}
\noindent \textbf{\boldmath{Example}} \textbf{\boldmath{2.}} \emph{Consider the variable fractional order force-free Duffing-Vander pol oscillator equation by replacing the constant fractional order~\cite{B26-mathematics-2384216} as}
\begin{equation}
\label{eq:FD432-mathematics-2384216}
\text{D}_{0,\text{t}}^{\text{$\upalpha$}(\text{t})}\Im(\text{t}) - \text{$\upmu$}\Im^{\prime}(\text{t}) + \text{$\upmu$}\Im^{\prime}(\text{t})\Im^{2}(\text{t}) + \text{a}\Im(\text{t}) + \text{b}\Im^{3}(\text{t}) = 0;\ \text{$\upalpha$}(\text{t}) \in \left( {1,2} \right\rbrack,
\tag{37}
\end{equation}
\emph{with the initial value conditions}
\begin{equation}
\nonumber\label{eq:FD433-mathematics-2384216}
\Im(0) = 2,\ \Im^{\prime}(0) = 0.
\end{equation}
\emph{We solve the example for} $\widetilde{\text{$\upsigma$}} = 4,6\ (\text{k} = 1,\ \text{M} = 3,5) $ \emph{by mentioned scheme and simulate the model for different parameters}.

\vspace{12pt}
In considering the problem, the following two cases of fractional order are considered:
\begin{enumerate}
\item[(i)] Constant order: $\text{$\upalpha$}(\text{t}) = 1.2,1.4,1.5,1.6,1.8$;
\item[(ii)] Variable order $\text{$\upalpha$}(\text{t}) = 1 + \sin \text{t}$.
\end{enumerate}

\begin{figure}[H]
\includegraphics[scale=1]{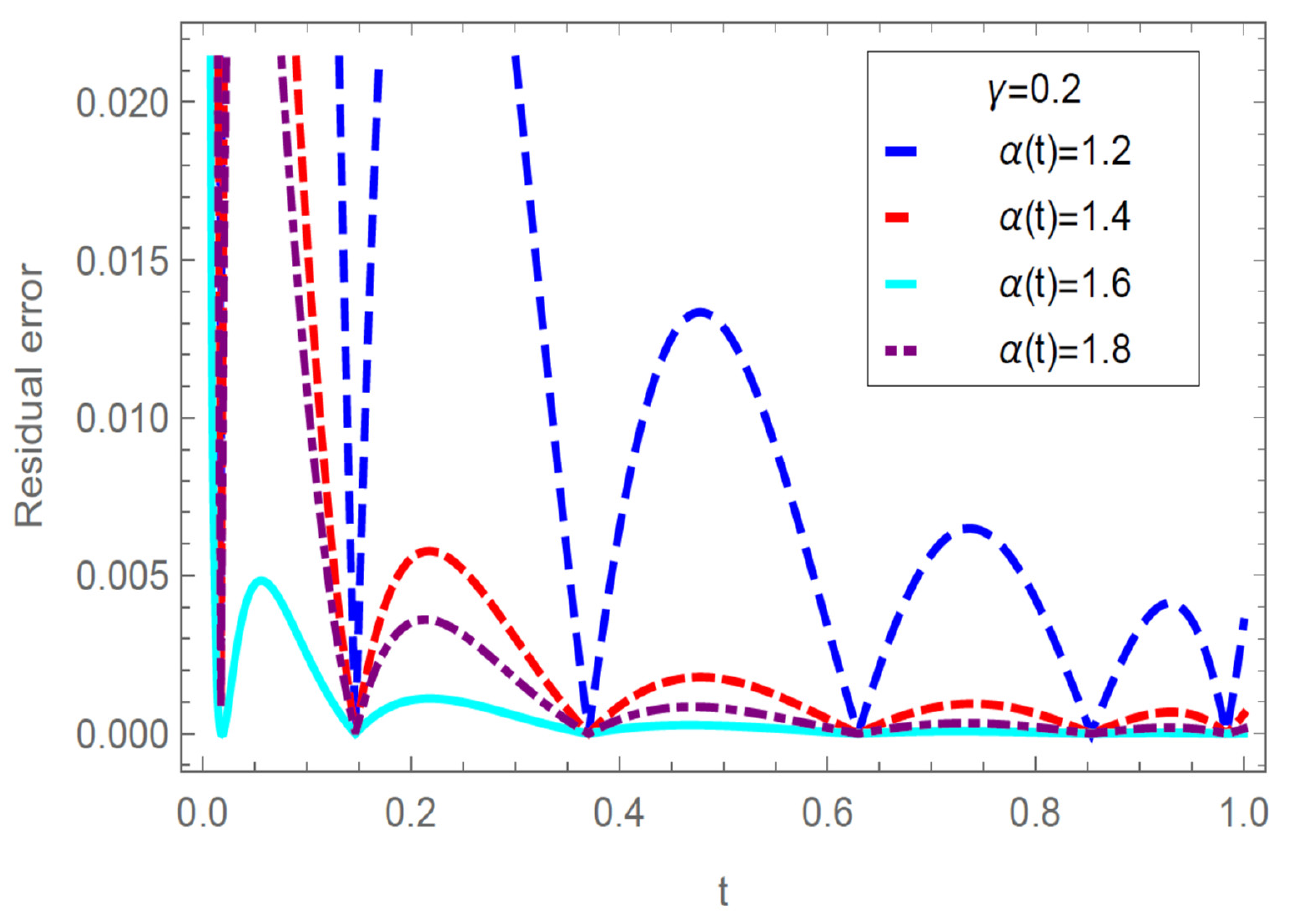}
\caption{Graph of Residual errors in Example 2 with $k = 1$, $M = 5$.}
\label{fig:mathematics-2384216-f006}
\end{figure}

The estimated AEs in the solutions of $\Im(\text{t}) $ with the comparison of adomian decomposition scheme (ADS) and restarted adomian decomposition scheme (RADS) for $\text{$\upalpha$}(\text{t}) = 2 $ and different FOBW bases are listed in Table~\ref{tabref:mathematics-2384216-t014}. It can be easily analyzed from Table~\ref{tabref:mathematics-2384216-t014} that the suitable value of $\text{$\upgamma$} $ is 1 for achieving the best accuracy in the solution of the given model, with $\text{$\upalpha$}(\text{t}) = 2 $ and the proposed approach is superior to ADS~\cite{B26-mathematics-2384216} and RADS~\cite{B27-mathematics-2384216} by considering the Lindsted scheme solution~\cite{B26-mathematics-2384216} as an approximated analytical solution. The residual errors in the solutions of $\Im(\text{t}) $ for $\text{$\upalpha$}(\text{t}) = 1.5 $ and $\text{$\upalpha$}(\text{t}) = 1 + \sin \text{t} $ under different parameters mentioned are presented 
in Tables~\ref{tabref:mathematics-2384216-t015} and \ref{tabref:mathematics-2384216-t016}, respectively. Furthermore, the estimated residual errors in the solutions for $\text{$\upgamma$} = 0.2 $ and selected $\text{$\upalpha$}(\text{t}) $ are given in Table~\ref{tabref:mathematics-2384216-t017}. The graphical interpretation of residual errors of solutions for selected values of $\text{$\upalpha$}(\text{t}) $, and $\text{$\upgamma$} = 0.2 $ is shown in Figure~\ref{fig:mathematics-2384216-f006}. The computed solutions are obtained for the first time with the variable order of the introduced model in terms of residual errors.    
\begin{table}[H]
\tablesize{\small}
\caption{Estimated AEs of solutions for $\text{$\upalpha$}(\text{t}) = 2 $ and selected $\text{$\upgamma$} $ in Example 2.}
\label{tabref:mathematics-2384216-t014}
\setlength{\cellWidtha}{\textwidth/6-2\tabcolsep-0in}
\setlength{\cellWidthb}{\textwidth/6-2\tabcolsep-0in}
\setlength{\cellWidthc}{\textwidth/6-2\tabcolsep-0in}
\setlength{\cellWidthd}{\textwidth/6-2\tabcolsep-0in}
\setlength{\cellWidthe}{\textwidth/6-2\tabcolsep-0in}
\setlength{\cellWidthf}{\textwidth/6-2\tabcolsep-0in}
        
\scalebox{1}[1]{\begin{tabularx}{\textwidth}{>{\centering\arraybackslash}m{\cellWidtha}>{\centering\arraybackslash}m{\cellWidthb}>{\centering\arraybackslash}m{\cellWidthc}>{\centering\arraybackslash}m{\cellWidthd}>{\centering\arraybackslash}m{\cellWidthe}>{\centering\arraybackslash}m{\cellWidthf}}
\toprule
\multicolumn{6}{>{\centering\arraybackslash}m{\cellWidtha + \cellWidthb + \cellWidthc + \cellWidthd + \cellWidthe + \cellWidthf+10\tabcolsep}}{\textbf{$\textbf{a}\ \textbf{=}\ \textbf{1}\textbf{,}\ \textbf{b}\ \textbf{=}\ \textbf{0.01}\textbf{,}\ \bm{\upmu}\ \textbf{=}\ \textbf{0.1} $}}\\
\cmidrule{1-6}
\multirow{2}{*}{\parbox{\cellWidtha}{\centering \textbf{T}}\vspace{-4pt}} & \multicolumn{3}{>{\centering\arraybackslash}m{\cellWidthb + \cellWidthc + \cellWidthd+4\tabcolsep}}{\textbf{
Proposed Approach, $\textbf{k}\ \textbf{=}\ \textbf{1}\textbf{,}\ \textbf{M}\ \textbf{=}\ \textbf{5} $}} & \multicolumn{2}{>{\centering\arraybackslash}m{\cellWidthe + \cellWidthf+2\tabcolsep}}{\textbf{
Reference Approach}}\\
\cmidrule{2-6}
& \textbf{$\bm{\upgamma}\ \textbf{=}\ \textbf{0.5} $} & \textbf{$\bm{\upgamma}\ \textbf{=}\ \textbf{0.9} $} & \textbf{$\bm{\upgamma}\ \textbf{=}\ \textbf{1} $} & \textbf{ ADS {\cite{B26-mathematics-2384216}}
} & \textbf{ RADS {\cite{B27-mathematics-2384216}}
}\\
\cmidrule{1-6}
0.1 & $3.3 \times 10^{- 6} $ & $2.1 \times 10^{- 6} $ & $2.3 \times 10^{- 6} $ & $2.4 \times 10^{- 3} $ & $2.4 \times 10^{- 3} $\\
\cmidrule{1-6}
0.3 & $2.0 \times 10^{- 5} $ & $4.4 \times 10^{- 6} $ & $3.8 \times 10^{- 6} $ & $2.2 \times 10^{- 3} $ & $2.2 \times 10^{- 3} $\\
\cmidrule{1-6}
0.5 & $2.1 \times 10^{- 5} $ & $2.1 \times 10^{- 6} $ & $8.4 \times 10^{- 7} $ & $1.5 \times 10^{- 3} $ & $1.5 \times 10^{- 3} $\\
\cmidrule{1-6}
0.7 & $3.0 \times 10^{- 5} $ & $6.9 \times 10^{- 7} $ & $2.1 \times 10^{- 6} $ & $6.2 \times 10^{- 4} $ & $2.2 \times 10^{- 4} $\\
\cmidrule{1-6}
0.9 & $9.2 \times 10^{- 5} $ & $5.3 \times 10^{- 5} $ & $5.1 \times 10^{- 5} $ & $1.4 \times 10^{- 3} $ & $1.3 \times 10^{- 3} $\\
\bottomrule
\end{tabularx}}
\end{table}
\vspace{-12pt}
    
\begin{table}[H]
\tablesize{\small}
\caption{Estimated residual errors of solutions for $\text{$\upalpha$}(\text{t}) = 1.5, $ $\text{k} = 1,\text{M} = 5, $ and selected $\text{$\upgamma$} $ in Example 2.}
\label{tabref:mathematics-2384216-t015}
\setlength{\cellWidtha}{\textwidth/7-2\tabcolsep-0in}
\setlength{\cellWidthb}{\textwidth/7-2\tabcolsep-0in}
\setlength{\cellWidthc}{\textwidth/7-2\tabcolsep-0in}
\setlength{\cellWidthd}{\textwidth/7-2\tabcolsep-0in}
\setlength{\cellWidthe}{\textwidth/7-2\tabcolsep-0in}
\setlength{\cellWidthf}{\textwidth/7-2\tabcolsep-0in}
\setlength{\cellWidthg}{\textwidth/7-2\tabcolsep-0in}
        
\scalebox{1}[1]{\begin{tabularx}{\textwidth}{>{\centering\arraybackslash}m{\cellWidtha}>{\centering\arraybackslash}m{\cellWidthb}>{\centering\arraybackslash}m{\cellWidthc}>{\centering\arraybackslash}m{\cellWidthd}>{\centering\arraybackslash}m{\cellWidthe}>{\centering\arraybackslash}m{\cellWidthf}>{\centering\arraybackslash}m{\cellWidthg}}
\toprule
\multicolumn{7}{>{\centering\arraybackslash}m{\cellWidtha + \cellWidthb + \cellWidthc + \cellWidthd + \cellWidthe + \cellWidthf + \cellWidthg+12\tabcolsep}}{\textbf{$\textbf{a}\ \textbf{=}\ \textbf{b}\ \textbf{=}\ \textbf{0.5}\textbf{,}\ \textbf{f}\ \textbf{=}\ \textbf{0.5}\textbf{,}\ \bm{\upmu}\ \textbf{=}\ \textbf{0.1}\textbf{,}\ \bm{\upomega}\ \textbf{=}\ \textbf{0.79} $}}\\
\cmidrule{1-7}
\textbf{T} & \textbf{$\bm{\upgamma}\ \textbf{=}\ \textbf{0.1} $} & \textbf{$\bm{\upgamma}\ \textbf{=}\ \textbf{0.2} $} & \textbf{$\bm{\upgamma}\ \textbf{=}\ \textbf{0.3} $} & \textbf{$\bm{\upgamma}\ \textbf{=}\ \textbf{0.5} $} & \textbf{$\bm{\upgamma}\ \textbf{=}\ \textbf{0.9} $} & \textbf{$\bm{\upgamma}\ \textbf{=}\ \textbf{1.0} $}\\
\cmidrule{1-7}
0.1 & $8.9 \times 10^{- 3} $ & $4.4 \times 10^{- 3} $ & $4.7 \times 10^{- 2} $ & $6.8 \times 10^{- 2} $ & $3.0 \times 10^{- 1} $ & $3.8 \times 10^{- 1} $\\
\cmidrule{1-7}
0.3 & $1.1 \times 10^{- 3} $ & $6.4 \times 10^{- 4} $ & $4.7 \times 10^{- 2} $ & $1.7 \times 10^{- 2} $ & $1.6 \times 10^{- 1} $ & $2.5 \times 10^{- 1} $\\
\cmidrule{1-7}
0.5 & $3.3 \times 10^{- 4} $ & $3.6 \times 10^{- 4} $ & $2.1 \times 10^{- 2} $ & $9.9 \times 10^{- 3} $ & $1.5 \times 10^{- 1} $ & $2.5 \times 10^{- 1} $\\
\cmidrule{1-7}
0.7 & $3.6 \times 10^{- 5} $ & $1.9 \times 10^{- 4} $ & $8.3 \times 10^{- 3} $ & $4.3 \times 10^{- 3} $ & $9.1 \times 10^{- 2} $ & $1.7 \times 10^{- 1} $\\
\cmidrule{1-7}
0.9 & $3.7 \times 10^{- 5} $ & $1.4 \times 10^{- 4} $ & $2.1 \times 10^{- 3} $ & $2.6 \times 10^{- 3} $ & $7.2 \times 10^{- 2} $ & $1.4 \times 10^{- 1} $\\
\bottomrule
\end{tabularx}}
\end{table}
\vspace{-12pt}
    
\begin{table}[H]
\tablesize{\small}
\caption{Estimated residual errors for $\text{$\upalpha$}(\text{t}) = 1 + \sin \text{t}, $ $\text{k} = 1,\text{M} = 5, $ and different values of $\text{$\upgamma$} $ in Example 2.}
\label{tabref:mathematics-2384216-t016}
\setlength{\cellWidtha}{\textwidth/7-2\tabcolsep-0in}
\setlength{\cellWidthb}{\textwidth/7-2\tabcolsep-0in}
\setlength{\cellWidthc}{\textwidth/7-2\tabcolsep-0in}
\setlength{\cellWidthd}{\textwidth/7-2\tabcolsep-0in}
\setlength{\cellWidthe}{\textwidth/7-2\tabcolsep-0in}
\setlength{\cellWidthf}{\textwidth/7-2\tabcolsep-0in}
\setlength{\cellWidthg}{\textwidth/7-2\tabcolsep-0in}
        
\scalebox{1}[1]{\begin{tabularx}{\textwidth}{>{\centering\arraybackslash}m{\cellWidtha}>{\centering\arraybackslash}m{\cellWidthb}>{\centering\arraybackslash}m{\cellWidthc}>{\centering\arraybackslash}m{\cellWidthd}>{\centering\arraybackslash}m{\cellWidthe}>{\centering\arraybackslash}m{\cellWidthf}>{\centering\arraybackslash}m{\cellWidthg}}
\toprule
\multicolumn{7}{>{\centering\arraybackslash}m{\cellWidtha + \cellWidthb + \cellWidthc + \cellWidthd + \cellWidthe + \cellWidthf + \cellWidthg+12\tabcolsep}}{\textbf{$\textbf{a}\ \textbf{=}\ \textbf{b}\ \textbf{=}\ \textbf{0.5}\textbf{,}\ \textbf{f}\ \textbf{=}\ \textbf{0.5}\textbf{,}\ \bm{\upmu}\ \textbf{=}\ \textbf{0.1}\textbf{,}\ \bm{\upomega}\ \textbf{=}\ \textbf{0.79} $}}\\
\cmidrule{1-7}
\textbf{T} & \textbf{$\bm{\upgamma}\ \textbf{=}\ \textbf{0.1} $} & \textbf{$\bm{\upgamma}\ \textbf{=}\ \textbf{0.2} $} & \textbf{$\bm{\upgamma}\ \textbf{=}\ \textbf{0.3} $} & \textbf{$\bm{\upgamma}\ \textbf{=}\ \textbf{0.5} $} & \textbf{$\bm{\upgamma}\ \textbf{=}\ \textbf{0.9} $} & \textbf{$\bm{\upgamma}\ \textbf{=}\ \textbf{1.0} $}\\
\cmidrule{1-7}
0.1 & $2.2 \times 10^{- 1} $ & $2.1 \times 10^{- 1} $ & $2.5 \times 10^{- 1} $ & $5.2 \times 10^{- 1} $ & $1.6 \times 10^{+ 0} $ & $1.9 \times 10^{+ 0} $\\
\cmidrule{1-7}
0.3 & $8.3 \times 10^{- 2} $ & $8.2 \times 10^{- 2} $ & $9.6 \times 10^{- 2} $ & $2.4 \times 10^{- 1} $ & $1.5 \times 10^{+ 0} $ & $2.2 \times 10^{+ 0} $\\
\cmidrule{1-7}
0.5 & $5.5 \times 10^{- 2} $ & $5.5 \times 10^{- 2} $ & $6.4 \times 10^{- 2} $ & $1.7 \times 10^{- 1} $ & $1.7 \times 10^{+ 0} $ & $2.6 \times 10^{+ 0} $\\
\cmidrule{1-7}
0.7 & $2.5 \times 10^{- 2} $ & $2.5 \times 10^{- 2} $ & $2.8 \times 10^{- 2} $ & $7.7 \times 10^{- 2} $ & $1.1 \times 10^{+ 0} $ & $1.8 \times 10^{+ 0} $\\
\cmidrule{1-7}
0.9 & $1.4 \times 10^{- 2} $ & $1.4 \times 10^{- 2} $ & $1.6 \times 10^{- 2} $ & $4.5 \times 10^{- 2} $ & $8.5 \times 10^{- 1} $ & $1.5 \times 10^{+ 0} $\\
\bottomrule
\end{tabularx}}
\end{table}
\vspace{-12pt}
    
\begin{table}[H]
\tablesize{\small}
\caption{Estimated residual errors of solutions for different $\text{$\upalpha$}(\text{t}) $, k = 1 and $\text{$\upgamma$} = 0.2 $ in Example 2.}
\label{tabref:mathematics-2384216-t017}
    
\begin{adjustwidth}{-\extralength}{0cm}
\setlength{\cellWidtha}{\fulllength/9-2\tabcolsep-0in}
\setlength{\cellWidthb}{\fulllength/9-2\tabcolsep-0in}
\setlength{\cellWidthc}{\fulllength/9-2\tabcolsep-0in}
\setlength{\cellWidthd}{\fulllength/9-2\tabcolsep-0in}
\setlength{\cellWidthe}{\fulllength/9-2\tabcolsep-0in}
\setlength{\cellWidthf}{\fulllength/9-2\tabcolsep-0in}
\setlength{\cellWidthg}{\fulllength/9-2\tabcolsep-0in}
\setlength{\cellWidthh}{\fulllength/9-2\tabcolsep-0in}
\setlength{\cellWidthi}{\fulllength/9-2\tabcolsep-0in}
        
\scalebox{1}[1]{\begin{tabularx}{\fulllength}{>{\centering\arraybackslash}m{\cellWidtha}>{\centering\arraybackslash}m{\cellWidthb}>{\centering\arraybackslash}m{\cellWidthc}>{\centering\arraybackslash}m{\cellWidthd}>{\centering\arraybackslash}m{\cellWidthe}>{\centering\arraybackslash}m{\cellWidthf}>{\centering\arraybackslash}m{\cellWidthg}>{\centering\arraybackslash}m{\cellWidthh}>{\centering\arraybackslash}m{\cellWidthi}}
\toprule
\multicolumn{9}{>{\centering\arraybackslash}m{\cellWidtha + \cellWidthb + \cellWidthc + \cellWidthd + \cellWidthe + \cellWidthf + \cellWidthg + \cellWidthh + \cellWidthi+16\tabcolsep}}{\textbf{
Residual Error, $\bm{\upgamma}\ \textbf{=}\ \textbf{0.2} $}}\\
\cmidrule{1-9}
\multirow{2}{*}{\parbox{\cellWidtha}{\centering \textbf{t}}\vspace{-4pt}} & \multicolumn{2}{>{\centering\arraybackslash}m{\cellWidthb + \cellWidthc+2\tabcolsep}}{\textbf{$\bm{\upalpha}(\textbf{t})\ \textbf{=}\ \textbf{1.2} $}} & \multicolumn{2}{>{\centering\arraybackslash}m{\cellWidthd + \cellWidthe+2\tabcolsep}}{\textbf{$\bm{\upalpha}(\textbf{t})\ \textbf{=}\ \textbf{1.4} $}} & \multicolumn{2}{>{\centering\arraybackslash}m{\cellWidthf + \cellWidthg+2\tabcolsep}}{\textbf{$\bm{\upalpha}(\textbf{t})\ \textbf{=}\ \textbf{1.6} $}} & \multicolumn{2}{>{\centering\arraybackslash}m{\cellWidthh + \cellWidthi+2\tabcolsep}}{\textbf{$\bm{\upalpha}(\textbf{t})\ \textbf{=}\ \textbf{1.8} $}}\\
\cmidrule{2-9}
 & \textbf{M = 3} & \textbf{M = 5} & \textbf{M = 3} & \textbf{M = 5} & \textbf{M = 3} & \textbf{M = 5} & \textbf{M = 3} & \textbf{M = 5}\\
\cmidrule{1-9}
0.1 & $7.7 \times 10^{- 1} $ & $8.0 \times 10^{- 2} $ & $7.7 \times 10^{- 1} $ & $1.5 \times 10^{- 2} $ & $2.2 \times 10^{+ 0} $ & $2.4 \times 10^{- 3} $ & $4.1 \times 10^{+ 0} $ & $1.0 \times 10^{- 2} $\\
\cmidrule{1-9}
0.3 & $6.0 \times 10^{- 2} $ & $2.1 \times 10^{- 2} $ & $3.2 \times 10^{- 2} $ & $3.0 \times 10^{- 3} $ & $5.6 \times 10^{- 2} $ & $5.6 \times 10^{- 4} $ & $7.8 \times 10^{- 2} $ & $1.7 \times 10^{- 3} $\\
\cmidrule{1-9}
0.5 & $7.8 \times 10^{- 2} $ & $1.2 \times 10^{- 2} $ & $2.6 \times 10^{- 1} $ & $1.7 \times 10^{- 3} $ & $3.5 \times 10^{- 1} $ & $2.5 \times 10^{- 4} $ & $4.3 \times 10^{- 1} $ & $7.9 \times 10^{- 4} $\\
\cmidrule{1-9}
0.7 & $8.9 \times 10^{- 3} $ & $5.6 \times 10^{- 3} $ & $1.3 \times 10^{- 2} $ & $8.0 \times 10^{- 4} $ & $1.5 \times 10^{- 2} $ & $7.5 \times 10^{- 5} $ & $1.6 \times 10^{- 2} $ & $2.9 \times 10^{- 4} $\\
\cmidrule{1-9}
0.9 & $1.0 \times 10^{- 1} $ & $3.3 \times 10^{- 3} $ & $1.0 \times 10^{- 1} $ & $5.3 \times 10^{- 4} $ & $1.1 \times 10^{- 1} $ & $2.2 \times 10^{- 5} $ & $1.1 \times 10^{- 1} $ & $1.5 \times 10^{- 4} $\\
\bottomrule
\end{tabularx}}
\end{adjustwidth}
\end{table}
\vspace{-12pt}
    
\section{Conclusions and Future Work}
\label{sect:sec9-mathematics-2384216}

In this study, we investigated the dynamical behavior of the non-linear model of oscillation equations under the variable-order fractional Caputo differential operator using a numerical framework based on FOBWs. In the present work, a well-organized FOBWs scheme has been successfully utilized for solving variable-order fractional force-free and forced oscillation equations arising in several areas of engineering and applied science. The described method uses fewer FOBWs, which produces better results. Some impacts of fractional parameters under variable order have been represented by tables and graphs. The accuracy of the suggested scheme is essentially affected by the order of the FOBWs. The present numerical simulations are compared with the known literature for integer order, and they are in good agreement. In conclusion, fractional derivative operators under variable order can be applied as a powerful tool for analyzing the dynamical behavior of several real-world problems.

The results achieved through our approach emphasized that:
\begin{enumerate}[label=$\bullet$]
\item The approach achieves accurate solutions that are easy to implement.
\item The accuracy of the error can be controlled and reduced by increasing the wavelet bases in the approximate solution.
\item Different types of nonlinear variable-order fractional models can easily be handled through the mentioned scheme.
\item The stability analysis of the proposed scheme for the solutions of the variable-order fractional oscillation equation is an interesting problem for future study.
\end{enumerate}


\vspace{6pt}
\authorcontributions{Conceptualization, A.R., B.P.J., M.P. and D.F.M.T.; 
Methodology, A.R., B.P.J., M.P. and D.F.M.T.; Validation, A.R., B.P.J. and D.F.M.T.; 
Formal analysis, A.R., B.P.J. and D.F.M.T.; 
Writing---original draft preparation, A.R. and D.F.M.T.; 
Writing---review and editing, A.R., B.P.J. and D.F.M.T.; 
Visualization, A.R., B.P.J. and D.F.M.T.; Supervision, A.R., B.P.J. and D.F.M.T.; 
Project administration, B.P.J. and D.F.M.T.; Funding acquisition, D.F.M.T. 
All authors have read and agreed to the published version of the manuscript.}

\funding{The research was partially funded by the 
Portuguese Foundation for Science and Technology (FCT), project UIDB/04106/2020.}

\informedconsent{Not applicable.}

\dataavailability{No data were used to support this study.}

\conflictsofinterest{The authors declare no conflict of interest.}


\begin{adjustwidth}{-\extralength}{0cm}

\reftitle{References}

\end{adjustwidth}

\begin{adjustwidth}{-\extralength}{0cm}
\PublishersNote{}
\end{adjustwidth}



\begin{thebibliography}{999}
	
\bibitem{B1-mathematics-2384216}
Podlubny, I. \emph{Fractional Differential Equations}; Academic Press: San Diego, CA, USA, 1999.

\bibitem{B2-mathematics-2384216}
Mainardi, F. \emph{Fractional Calculus and Waves in Linear Viscoelasticity}; Imperial College: London, UK, 2010.

\bibitem{B3-mathematics-2384216}
Rossikhin, Y.A.; Shitikova, M.V. Application of fractional calculus for dynamic problems of solid mechanics: Novel trends and recent results. \emph{Appl. Mech. Rev.} \textbf{\boldmath{2010}}, \emph{63}, 010801.

\bibitem{B4-mathematics-2384216}
Ampun, S.; Sawangtong, P. The Approximate Analytic Solution of the Time-Fractional Black-Scholes Equation with a European Option Based on the Katugampola Fractional Derivative. \emph{Mathematics} \textbf{\boldmath{2021}}, \emph{9}, 214. [\href{https://doi.org/10.3390/math9030214}{CrossRef}]

\bibitem{B5-mathematics-2384216}
Alshbool, M.H.T.; Bataineh, A.S.; Hashim, I.; RasitIsik, O. Solution of fractional-order differential equations based on the operational matrices of new fractional Bernstein functions. \emph{J. King Saud Univ.-Sci.} \textbf{\boldmath{2017}}, \emph{29}, 1--18. [\href{https://doi.org/10.1016/j.jksus.2015.11.004}{CrossRef}]

\bibitem{B6-mathematics-2384216}
Alshbool, M.H.T.; Mohammad, M.; Isik, O.; Hashim, I. Fractional Bernstein operational matrices for solving integro-differential equations involved by Caputo fractional derivative. \emph{Results Appl. Math.} \textbf{\boldmath{2022}}, \emph{14}, 100258. [\href{https://doi.org/10.1016/j.rinam.2022.100258}{CrossRef}]

\bibitem{B7-mathematics-2384216}
Rogosin, S.; Dubatovskaya, M. Letnikov vs. Marchaud: A survey on two prominent constructions of fractional derivatives. \emph{Mathematics} \textbf{\boldmath{2018}}, \emph{6}, 3. [\href{https://doi.org/10.3390/math6010003}{CrossRef}]

\bibitem{B8-mathematics-2384216}
Cai, M.; Li, C. Numerical Approaches to Fractional Integrals and Derivatives: A Review. \emph{Mathematics} \textbf{\boldmath{2020}}, \emph{8}, 43. [\href{https://doi.org/10.3390/math8010043}{CrossRef}]

\bibitem{B9-mathematics-2384216}
Sun, H.G.; Chen, W.; Sheng, H.; Chen, Y.Q. On mean square displacement behaviors of anomalous diffusions with variable and random orders. \emph{Phys. Lett. A} \textbf{\boldmath{2010}}, \emph{374}, 906--910. [\href{https://doi.org/10.1016/j.physleta.2009.12.021}{CrossRef}]

\bibitem{B10-mathematics-2384216}
Ingman, D.; Suzdalnitsky, J.; Zeifman, M. Constitutive dynamic order model for nonlinear contact phenomena. \emph{J. Appl. Mech.} \textbf{\boldmath{2000}}, \emph{67}, 383--390. [\href{https://doi.org/10.1115/1.1304916}{CrossRef}]

\bibitem{B11-mathematics-2384216}
Sun, H.G.; Zhang, H.; Chen, W.; Reeves, D.M. Use of a variable index fractional-derivative model to capture transient dispersion in heterogeneous media. \emph{J. Contam. Hydrol.} \textbf{\boldmath{2014}}, \emph{157}, 47--58. [\href{https://doi.org/10.1016/j.jconhyd.2013.11.002}{CrossRef}]

\bibitem{B12-mathematics-2384216}
Gomez-Aguilar, J.F. Analytical and Numerical solutions of a nonlinear alcoholism model via variable-order fractional differential equations. \emph{J. Phys. A} \textbf{\boldmath{2018}}, \emph{494}, 52--75. [\href{https://doi.org/10.1016/j.physa.2017.12.007}{CrossRef}]

\bibitem{B13-mathematics-2384216}
Debnath, L. \emph{Wavelets Transform and Their Applications}; Birkhauser: Boston, MA, USA, 2002.

\bibitem{B14-mathematics-2384216}
Chui, C.K. \emph{An Introduction to Wavelets}; Academic Press: San Diego, CA, USA, 1992.

\bibitem{B15-mathematics-2384216}
Rayal, A.; Verma, S.R. An approximate wavelets solution to the class of variational problems with fractional order. \emph{J. Appl. Math. Comput.} \textbf{\boldmath{2020}}, \emph{65}, 735--769. [\href{https://doi.org/10.1007/s12190-020-01413-9}{CrossRef}]

\bibitem{B16-mathematics-2384216}
Rayal, A.; Verma, S.R. Numerical study of variational problems of moving or fixed boundary conditions by Muntz wavelets. \emph{J. Vib. Control} \textbf{\boldmath{2020}}, \emph{28}, 214--229. [\href{https://doi.org/10.1177/1077546320974792}{CrossRef}]

\bibitem{B17-mathematics-2384216}
Rayal, A.; Verma, S.R. Numerical analysis of pantograph differential equation of the stretched type associated with fractal-fractional derivatives via fractional order Legendre wavelets. \emph{Chaos Solitons Fractals} \textbf{\boldmath{2020}}, \emph{139}, 110076. [\href{https://doi.org/10.1016/j.chaos.2020.110076}{CrossRef}]

\bibitem{B18-mathematics-2384216}
Rayal, A.; Verma, S.R. Two-dimensional Gegenbauer wavelets for the numerical solution of tempered fractional model of the nonlinear Klein-Gordon equation. \emph{Appl. Numer. Math.} \textbf{\boldmath{2022}}, \emph{174}, 191--220. [\href{https://doi.org/10.1016/j.apnum.2022.01.015}{CrossRef}]

\bibitem{B19-mathematics-2384216}
Rayal, A.; Tamta, S.; Rawat, S.; Kashif, M. Numerical view of Lucas-Lehmer polynomials with its characteristics. \emph{Uttaranchal J. Appl. Life Sci. Uttaranchal Univ.} \textbf{\boldmath{2022}}, \emph{3}, 66--75.

\bibitem{B20-mathematics-2384216}
Rayal, A. An effective Taylor wavelets basis for the evaluation of numerical differentiations. \emph{Palest. J. Math.} \textbf{\boldmath{2023}}, \emph{12}, 551--568.

\bibitem{B21-mathematics-2384216}
\textls[-45]{Almeida, R.; Tavares, D.; Torres, D.F.M. \emph{The Variable-Order Fractional Calculus of Variations}; Springer: Cham, Switzerland, 2019.} [\href{https://doi.org/10.1007/978-3-319-94006-9}{CrossRef}]

\bibitem{B22-mathematics-2384216}
Abramowitz, M.; Stegun, I. \emph{Handbook of Mathematical Functions}; Dover: New York, NY, USA, 1964.

\bibitem{B23-mathematics-2384216}
Chua, L.O.; Desoer, C.A.; Kuh, E.S. \emph{Linear and Nonlinear Circuits, McGraw-Hill Series in Electrical Engineering: Circuits and Systems}; McGraw-Hill: New York, NY, USA, 1987.

\bibitem{B24-mathematics-2384216}
Barbosa, R.S.; Machado, J.A.T.; Vinagre, B.M.; Calderon, A.J. Analysis of the Van der Pol oscillator containing derivatives of fractional order. \emph{J. Vib. Control} \textbf{\boldmath{2007}}, \emph{1}, 1291--1301. [\href{https://doi.org/10.1177/1077546307077463}{CrossRef}]

\bibitem{B25-mathematics-2384216}
Kyamakya, K.; Ngoy, C.; Tamasala, M.; Chedjou, J. A novel image processing approach combining a ‘coupled nonlinear oscillators’-based paradigm with cellular neural networks for dynamic robust contrast enhancement. In \emph{ISAST Transactions on Computers and Intelligent Systems, Proceedings of the 12th International Workshop on Cellular Nanoscale Networks and Their Applications (CNNA 2010), Berkeley, CA, USA, 3--5 February 2010}; IEEE: New York, NY, USA, 2010; pp. 1--7.

\bibitem{B26-mathematics-2384216}
Cordshooli, G.A.; Vahidi, A.R. Solution of Duffing-van der pol equation using decomposition method. \emph{Adv. Stud. Theor. Phys.} \textbf{\boldmath{2011}}, \emph{5}, 121--129.

\bibitem{B27-mathematics-2384216}
Vahidi, A.R.; Azimzadeh, Z.; Mohammadifar, S. Restarted Adomian Decomposition Method for Solving Duffing-van der Pol Equation. \emph{Appl. Math. Sci.} \textbf{\boldmath{2012}}, \emph{6}, 499--507.

\bibitem{B28-mathematics-2384216}
Doha, E.H.; Abd-Elhameed, W.M.; Youssri, Y.H. New ultraspherical wavelets collocation method for solving 2nth-order initial and boundary value problems. \emph{J. Egypt. Math. Soc.} \textbf{\boldmath{2015}}, \emph{36}, 319--327. [\href{https://doi.org/10.1016/j.joems.2015.05.002}{CrossRef}]

\bibitem{B29-mathematics-2384216}
Mohyud-Din, S.T.; Iqbal, M.A.; Hassan, S.M. Modified Legendre Wavelets Technique for Fractional Oscillation Equations. \emph{Entropy} \textbf{\boldmath{2015}}, \emph{17}, 6925--6936. [\href{https://doi.org/10.3390/e17106925}{CrossRef}]

\bibitem{B30-mathematics-2384216}
Khan, M.M.-U.-R. Analytical Solution of Van Der Pol’s Differential Equation Using Homotopy Perturbation Method. \emph{J. Appl. Math. Phys.} \textbf{\boldmath{2019}}, \emph{7}, 1--12. [\href{https://doi.org/10.4236/jamp.2019.71001}{CrossRef}]

\bibitem{B31-mathematics-2384216}
Kumar, M.; Varshney, P. Numerical Simulation of Van der Pol Equation Using Multiple Scales Modified Lindstedt--Poincare Method. \emph{Proc. Natl. Acad. Sci. India Sect. A Phys. Sci.} \textbf{\boldmath{2021}}, \emph{91}, 55--65. [\href{https://doi.org/10.1007/s40010-019-00655-y}{CrossRef}]

\bibitem{B32-mathematics-2384216}
Hamed, M.; El-Kalla, I.; El-Beltagy, M.; El-Desouky, B. Numerical solutions of stochastic Duffing-Van der Pol equations. \emph{Indian J. Pure Appl. Math.} \textbf{\boldmath{2023}}. [\href{https://doi.org/10.1007/s13226-022-00361-3}{CrossRef}]

\bibitem{B33-mathematics-2384216}
Bhrawy, A.H.; Zaky, M.A. Numerical simulation for two-dimensional variable-order fractional nonlinear cable equation. \emph{Nonlinear Dyn.} \textbf{\boldmath{2015}}, \emph{80}, 101--116. [\href{https://doi.org/10.1007/s11071-014-1854-7}{CrossRef}]

\bibitem{B34-mathematics-2384216}
Bhrawy, A.H.; Zaky, M.A.; Alzaidy, J.F. two shifted Jacobi-Gauss collocation schemes for solving two-dimensional variable-order fractional Rayleigh-Stokes problem. \emph{Adv. Diff. Equ.} \textbf{\boldmath{2016}}, \emph{2016}, 272. [\href{https://doi.org/10.1186/s13662-016-0998-9}{CrossRef}]

\bibitem{B35-mathematics-2384216}
Bhrawy, A.H.; Zaky, M.A. An improved collocation method for multi-dimensional space-time variable-order fractional Schrödinger equations. \emph{Appl. Numer. Math.} \textbf{\boldmath{2017}}, \emph{111}, 197--218. [\href{https://doi.org/10.1016/j.apnum.2016.09.009}{CrossRef}]

\bibitem{B36-mathematics-2384216}
Xu, Y.; Erturk, V.S. A finite difference technique for solving variable-order fractional integro-differential equations. \emph{Bull. Iran. Math. Soc.} \textbf{\boldmath{2014}}, \emph{40}, 699--712.

\bibitem{B37-mathematics-2384216}
Wang, Z.; Vong, S. Compact difference schemes for the modified anomalous fractional sub-diffusion equation and the fractional diffusion-wave equation. \emph{J. Comput. Phys.} \textbf{\boldmath{2014}}, \emph{277}, 1--15. [\href{https://doi.org/10.1016/j.jcp.2014.08.012}{CrossRef}]

\bibitem{B38-mathematics-2384216}
Fu, Z.J.; Chen, W.; Ling, L. Method of approximate particular solutions for constant and variable-order fractional diffusion models. \emph{Eng. Anal. Bound. Elem.} \textbf{\boldmath{2015}}, \emph{57}, 37--46. [\href{https://doi.org/10.1016/j.enganabound.2014.09.003}{CrossRef}]

\bibitem{B39-mathematics-2384216}
Zayernouri, M.; Karniadakis, G.E. Fractional spectral collocation methods for linear and nonlinear variable order FPDEs. \emph{J.~Comput. Phys.} \textbf{\boldmath{2015}}, \emph{293}, 312--338. [\href{https://doi.org/10.1016/j.jcp.2014.12.001}{CrossRef}]

\bibitem{B40-mathematics-2384216}
Chen, Y.M.; Wei, Y.Q.; Liu, D.Y.; Yu, H. Numerical solution for a class of nonlinear variable order fractional differential equations with Legendre wavelets. \emph{Appl. Math. Lett.} \textbf{\boldmath{2015}}, \emph{46}, 83--88. [\href{https://doi.org/10.1016/j.aml.2015.02.010}{CrossRef}]

\bibitem{B41-mathematics-2384216}
Yaghoobi, S.; Moghaddam, B.P.; Ivaz, K. An efficient cubic spline approximation for variable-order fractional differential equations with time delay. \emph{Nonlinear Dyn.} \textbf{\boldmath{2017}}, \emph{87}, 815--826. [\href{https://doi.org/10.1007/s11071-016-3079-4}{CrossRef}]

\bibitem{B42-mathematics-2384216}
Zuniga-Aguilar, C.J.; Coronel-Escamilla, A.; Gomez-Aguilar, J.F.; Alvarado-Martinez, V.M.; Romero-Ugalde, H.M. New numerical approximation for solving fractional delay differential equations of variable order using artificial neural networks. \emph{Eur. Phys. J. Plus} \textbf{\boldmath{2018}}, \emph{133}, 75. [\href{https://doi.org/10.1140/epjp/i2018-11917-0}{CrossRef}]

\bibitem{B43-mathematics-2384216}
Heydari, M.H. Chebyshev cardinal functions for a new class of nonlinear optimal control problems generated by atangana baleanu caputo variable-order fractional derivative. \emph{Chaos Solitons Fractal} \textbf{\boldmath{2020}}, \emph{130}, 109401. [\href{https://doi.org/10.1016/j.chaos.2019.109401}{CrossRef}]

\bibitem{B44-mathematics-2384216}
Nemati, S.; Lima, P.M.; Torres, D.F.M. Numerical Solution of Variable-Order Fractional Differential Equations Using Bernoulli Polynomials. \emph{Fractal Fract.} \textbf{\boldmath{2021}}, \emph{5}, 219. [\href{https://doi.org/10.3390/fractalfract5040219}{CrossRef}]

\bibitem{B45-mathematics-2384216}
Kaabar, M.K.A.; Refice, A.; Souid, M.S.; Mart{\fontencoding{T5}\selectfont{\'i}}nez, F.; Etemad, S.; Siri, Z.; Rezapour, S. Existence and U-H-R Stability of Solutions to the Implicit Nonlinear FBVP in the Variable Order Settings. \emph{Mathematics} \textbf{\boldmath{2021}}, \emph{9}, 1693. [\href{https://doi.org/10.3390/math9141693}{CrossRef}]

\bibitem{B46-mathematics-2384216}
Mirzaee, F.; Alipour, S. Fractional-order orthogonal Bernstein polynomials for numerical solution of nonlinear fractional partial Volterra integro-differential equations. \emph{Math. Meth. Appl. Sci.} \textbf{\boldmath{2019}}, \emph{42}, 1870--1893. [\href{https://doi.org/10.1002/mma.5481}{CrossRef}]

\bibitem{B47-mathematics-2384216}
Wang, J.S.; Liu, L.Q.; Liu, L.C.; Chen, Y.M. Numerical solution for the variable order fractional partial differential equation with Bernstein polynomials. \emph{Int. J. Adv. Comput. Technol.} \textbf{\boldmath{2014}}, \emph{6}, 22--37.

\end{thebibliography}
\end{document}